\documentclass[10pt, reqno]{amsart}

\usepackage{setspace} 
\usepackage{amssymb}
\usepackage{amsmath}
\usepackage{amsthm}
\usepackage[a4paper]{geometry}
\usepackage{amsfonts}
\usepackage{bbm}
\usepackage{xcolor}
\newcommand{\brac}[1]{\left( #1 \right)}
\usepackage{changes}

\usepackage{color}
\usepackage[colorlinks=true,linkcolor=blue,citecolor=blue,pdfpagelabels=false]{hyperref}
\usepackage[nameinlink]{cleveref}

\usepackage{url}

\makeatletter
\newcommand{\myitem}[2]{%
  \item[#1]%
  \protected@edef\@currentlabel{#1}%
  \label{#2}%
}
\makeatother

\allowdisplaybreaks[4] 
\usepackage{aliascnt}

\newtheorem{thm}{Theorem}[section]

\newaliascnt{lem}{thm}
\newtheorem{lem}[lem]{Lemma}
\aliascntresetthe{lem}

\newaliascnt{prop}{thm}
\newtheorem{prop}[prop]{Proposition}
\aliascntresetthe{prop}

\newaliascnt{cor}{thm}
\newtheorem{cor}[cor]{Corollary}
\aliascntresetthe{cor}

\theoremstyle{definition}
\newaliascnt{example}{thm}
\newtheorem{example}[example]{Example}
\aliascntresetthe{example}
\theoremstyle{definition}

\newaliascnt{defn}{thm}
\newtheorem{defn}[defn]{Definition}
\aliascntresetthe{defn}

\newaliascnt{rem}{thm}
\newtheorem{rem}[rem]{Remark}
\aliascntresetthe{rem}

\crefname{thm}{theorem}{theorems}
\Crefname{thm}{Theorem}{Theorems}

\crefname{lem}{lemma}{lemmas}
\Crefname{lem}{Lemma}{Lemmas}

\crefname{rem}{remark}{remarks}
\Crefname{rem}{Remark}{Remarks}

\crefname{example}{example}{examples}
\Crefname{example}{example}{examples}

\def\C{\mathbb C}

\def\N{\mathbb N}

\def\R{\mathbb R}

\newcommand{\CB}{\mathcal{B}}

\newcommand{\CD}{\mathcal{D}}

\newcommand{\CF}{\mathcal{F}}

\newcommand{\CH}{\mathcal{H}}

\newcommand{\CM}{\mathcal{M}}

\newcommand{\CP}{\mathcal{P}}

\newcommand{\auleq}{\underset{a.u.}{\leq }}
\newcommand{\bauleq}{\underset{b.a.u.}{\leq }}

\usepackage{enumitem}    

\def\be{\begin{equation}}
	\def\ee{\end{equation}}
\def\bt{\begin{Theorem}}
	\def\et{\end{Theorem}}
\def\bi{\begin{itemize}}
	\def\ei{\end{itemize}}
\def\bea{\begin{eqnarray}}
	\def\eea{\end{eqnarray}}
\def\beast{\begin{eqnarray*}}
	\def\eeast{\end{eqnarray*}}
\def\ben{\begin{enumerate}}
	\def\een{\end{enumerate}}

\newcommand{\norm}[1]{\left\Vert#1\right\Vert}

\newcommand{\lrp}[1]{\left(#1 \right)}
\newcommand{\Wstar}[1]{W^{*}\left(#1\right)}
\newcommand{\normi}[1]{\left\lVert #1 \right\rVert_{\infty}}

\usepackage{amsthm}
\theoremstyle{plain}
\newtheorem{mainthm}{Theorem}

\renewcommand{\MR}[1]{} 

\newcommand{\abs}[1]{\left\vert#1\right\vert}

\numberwithin{equation}{section}

\allowdisplaybreaks[4]

\makeatletter
\renewcommand{\@biblabel}[1]{\textbf{[#1]}}
\makeatother

\begin{document}

\title{Weighted Logarithm Laws in non-commutative Probability}

	\author[S. Panja]{Sourav Panja}
	\address{Theoretical Statistics and Mathematics Unit, Indian Statistical Institute, Delhi, India}
	\email{souravpanja07@gmail.com}
	
	\author[D. Saha]{Diptesh Saha}
	\address{Theoretical Statistics and Mathematics Unit, Indian Statistical Institute, Kolkata, India}
	\email{dptshs@gmail.com}

\markboth
{Sourav Panja and Diptesh Saha}
{Non-commutative analouge of Tomkins LIL for Martingale differences}



\begin{abstract}
We prove analogues of weighted laws of iterated logarithm and laws of logarithm for non-commutative martingales and independent random variables respectively. First we obtain a weighted LIL for non-commutative martingales, extending the sharp unweighted result of \cite{panja2026non} to the case for polynomial weights, then by approximation to a broader class of functions. Then we proved a weighted law of logarithm for triangular array of weights and replacing the coefficient in the denominator by regularly varying functions with index $\alpha$. We show that the normalised sums converge to 0 for $\alpha> 1/2$ and for the case $\alpha=1/2$, we obtain  the sharp upper bound with constant $\sqrt{2}$.

\end{abstract}

\maketitle

\section{Introduction}

The law of the iterated logarithm (LIL) is one of the most important and precise limit theorems in probability theory. The strong law of large numbers establishes the almost sure convergence of normalized sums and the central limit theorem describes their limiting distribution after rescaling. The LIL tells us how large these fluctuations can be infinitely often. To see the basic idea, let $S_n=\sum_{k=1}^n X_k,$  where $(X_k)$ is a sequence of centered independent random variables. In the classical finite-variance case, the LIL uses the scale $\sqrt{2s_n^2\log\log s_n^2},$ where $s_n^2=\sum_{k=1}^n \mathbb E(X_k^2).$ For independent identically distributed (i.i.d) random variables with variance one, this becomes $\sqrt{2n\log\log n}.$ The term $\log\log n$ is the main feature of the theorem.

The history of the LIL goes back to the early work of Khintchine and Kolmogorov \cite{kolmogoroff1929gesetz}. Kolmogorov proved an LIL for independent random variables under suitable conditions on the individual summands. Hartman and Wintner later obtained the classical finite-variance result for i.i.d random variables \cite{hartman1941law}. Their work showed that, in the i.i.d. setting, boundedness of the random variables is not needed. The LIL was then extended in several directions. Stout developed martingale versions of the LIL \cite{stout1970martingale}. This was an important step. Martingales are not independent, but they still have a useful conditional mean-zero structure. Hence, martingale LILs provide a natural way to study dependent random sequences.

Another important direction is the study of weighted sums. Instead of the ordinary partial sum $S_n=\sum_{k=1}^n X_k,$ one considers $S_n(f)=\sum_{k=1}^n f(k/n)X_k,$ where $f$ is a function on $[0,1]$. The case $f\equiv1$ gives the ordinary partial sum. Thus, weighted sums contain the classical LIL as a special case. 

Weighted LILs have a long history. Gaposhkin studied LILs connected with Ces\`aro and Abel summation methods \cite{gapo1965the}. Tomkins studied weighted averages of independent random variables \cite{tomkins1971iterated}. For example, for suitable functions $f$, Tomkins proved in \cite{tomkins1972a} that
$$\limsup_{n\to\infty} \frac{1}{\sqrt{2n\log\log n}} \sum_{k=1}^n f(k/n)X_k = \left(\int_0^1 f(t)^2\,dt\right)^{1/2} \qquad\text{a.s.}$$
when the random variables are centered and have variance one. He later extended this result to martingale difference sequences \cite{tomkins1975iterated}. Li and Tomkins later obtained further results for weighted sums of i.i.d. random variables and also studied Banach-space versions \cite{li1996laws}. Their work shows that the classical LIL also holds for a broad class of regular weights.

This leads to a natural question: \textit{What happens if the weights are allowed to be much more general?} This question is important because the weights in $f(k/n)$ have a special structure. A general triangular array does not have this property. Let $A=(a_{n,k}),~ 1\leq k\leq n,$ and consider $S_n(A)=\sum_{k=1}^n a_{n,k}X_k.$ The coefficients may now change in an arbitrary way from one row to the next. Thus, $S_n(A)$ can behave differently from ordinary partial sums. This difference is serious. The classical LIL need not hold for bounded triangular arrays. Li and Rao constructed a bounded triangular array for which 
$$\limsup_{n\to\infty} \frac{S_n(A)} {\sqrt{2n\log\log n}} =\infty \qquad\text{a.s.},$$ showing that the classical LIL fails \cite{deli1996note}. Thus, even bounded weights can produce fluctuations much larger than the classical LIL scale. This shows that boundedness of the weights alone is not enough to preserve the $\log\log n$ scale.

This observation changes the question. For general triangular arrays, we must first ask whether the LIL is even the correct scale. Li, Rao, and Wang first studied this problem for weighted sums with multidimensional indices \cite{deli1995strong}. In the one-dimensional case, they assumed $\mathbb E(X^2)<\infty,$ and considered a bounded triangular array of weights $K:=\sup_{n\geq1}\max_{1\leq k\leq n}|a_{n,k}|<\infty.$ They showed that the natural scale is $\sqrt{2n\log n}$ and proved 
\begin{equation}\label{eqn:Li et al law of logarithm}
    \limsup_{n\to\infty} \frac{|S_n(A)|} {\sqrt{2n\log n}} \leq 12K\lrp{\mathbb E X^2}^{1/2} \qquad\text{a.s.}
\end{equation}
They call \eqref{eqn:Li et al law of logarithm} the law of logarithms. Later, Chen and Hao \cite{chen2011remark} obtained the sharper upper bound $K\sqrt{\mathbb E X^2}$ under the same finite second-moment condition. 

Li and Tomkins \cite{li2003law} generalize \eqref{eqn:Li et al law of logarithm}. They replaced the specific $\sqrt{2n\log n}$ normalization by more general $b_n=B(n)$, where $B$ is regularly varying with index $\alpha\in[1/2,1)$, without assuming $\mathbb E(X^2)<\infty$. When $\alpha>1/2$, they proved that the normalized weighted sums converge to zero, while at the critical value $\alpha=1/2$ they obtained a non-trivial upper bound:
\begin{equation}\label{eqn:Li tomkins LL bound}
    \limsup_{n\to\infty}\frac{|S_n(A)|}{b_n} \leq 6\sqrt{2}\,\ell\sup_{1\leq k\leq n}|a_{n,k}| \qquad\text{ a.s.}
\end{equation}
under a suitable truncated second-moment condition. Its non-commutative counterpart is formulated in Condition~\eqref{item:alpha half condition} of \Cref{mainthm:LIL deterministic weights}.

Over the past several decades, non-commutative probability theory has developed into a rich generalization of the classical theory. In this framework, the underlying probability space is replaced by a finite von Neumann algebra $(\CM,\tau)$ with a faithful normal tracial state, and random variables by $\tau$-measurable operators. The non-commutative setting is not a simple change of notation. Operators do not commute. Hence many pointwise arguments from classical probability are no longer available.

Non-commutative martingale convergence and ergodic theorems began with the original work of Cuculescu \cite{cuculescu1971martingales}. Junge and Xu \cite{junge2002doob, JX07} further developed these concepts by focusing on martingale and maximal inequalities. However, progress on the non-commutative LIL has been more recent. Konwerska \cite{konwerska2008law} obtained a result similar to the Hartman-Wintner theorem for independent self-adjoint operators. Zeng proved a non-commutative analogue of Stout's martingale LIL \cite{zeng2015kolmogorov}. His result gave an upper bound with constant $2$. A sharp version was obtained in our companion paper \cite{panja2026non}. There we proved optimal non-commutative analogues of the classical LIL for martingales and for independent random variables.

The goal of this paper is to develop the weighted versions of these results in a non-commutative framework. Here, $(\CM,\tau)$ represents a finite von Neumann algebra with an increasing sequence $(\CM_n)$ of von Neumann subalgebras. For $t > e$, set $L(t) := \max\{1, \log\log t\}$. Our first main result presents the following non-commutative version of Tomkins' theorem \cite{tomkins1975iterated}. 
\begin{mainthm}\label{mainthm: polynomial bound tomkins}
    Let $(\CM, \tau)$ be a finite von Neumann algebra and $(d_n)$ be a sequence of self-adjoint martingale differences  with $\tau(d_1)=0$, $E_{n-1}(d_n^2)=1$ and $\max_{m \le n} \norm{d_m} \le M_n$, where $(M_n)$ is a sequence of positive reals such that $\lim\limits_{n \to \infty} \frac{L(n) M_n^2}{n}=0.$
	Let $p$ be a polynomial on $[0,1]$, then
    \begin{align}\label{eqn: polynomial eqn}
        \limsup_{n \to \infty} \frac{1}{(\sqrt{2n L(n)})} \sum_{m=1}^n p(m/n) d_m \underset{a.u.}{\le } \|p\|_2.
    \end{align}
\end{mainthm}
When $p\equiv1$, this gives the non-commutative LIL from \cite{panja2026non}. Thus, the present theorem extends the earlier unweighted result. The first difficulty is the triangular form of the weights. The coefficient $p(m/n)$ depends on $n$. Therefore, the weighted sum is not an ordinary martingale sum. A direct use of the unweighted LIL is not enough. In the classical theory, Abel summation is a useful tool for weighted sums. In the non-commutative setting, however, the same argument cannot be used directly. We need estimates that are compatible with almost uniform convergence and with the lack of commutativity.

Our proof starts with monomial weights $p(t)=t^j. $ The corresponding sums are $\sum_{m=1}^n m^j d_m.$  We then combine these estimates for finitely many monomials to treat polynomial weights. A blocking argument is used throughout this step. We divide the indices into geometric blocks. On one block, the relevant normalizing quantities are comparable. This allows the triangular weights to be compared with fixed coefficients. We then use non-commutative maximal estimates and a non-commutative Borel-Cantelli argument. After the polynomial case, we pass to a suitable class of functions of bounded variation by approximation. We also obtain the result for the power series considered in the paper. Thus, the polynomial case is the basic step, while approximation gives the larger class of weights.

The second part of the paper asks a different question: \textit{What happens for general bounded deterministic weights?} This is the natural non-commutative counterpart of the class studied by Li and Tomkins \cite{li2003law}.
\begin{mainthm}\label{mainthm:LIL deterministic weights}
    Let $(x_k)_{k\ge 1}$ be a i.i.d sequence of mean zero self-adjoint random variables on $(\CM,\tau)$. Let $A=(a_{n,k})\in\CD$ with $K=\sup_{n,k}|a_{n,k}|$, and let $B$, $(b_n)$ be as above with index $\alpha\in[\tfrac12,1)$. Assume that,
    \begin{enumerate}
        \item\label{item:first condition} $\tau\lrp{B^{-1}(|x_1|)} < \infty$
        \item\label{item:second condition} $\frac{B^{-1}(t\log t)}{\log t} G(t) \;\longrightarrow\; 0$ as $t\to\infty$,
        \item\label{item:alpha half condition} if $\alpha= \tfrac{1}{2}$, then there exists a constant $0\le \ell <\infty$ such that
              $$\limsup_{n\to\infty} \; \frac{n\log n}{b_n^2} H\left(\frac{b_n}{\log n}\right) = \ell^2.$$
    \end{enumerate}
    Then:
    \begin{enumerate}
        \item[(a)] If $\tfrac12<\alpha<1$ and $x_1\in L_2(\CM,\tau)$, $\displaystyle \frac{S_n(A)}{b_n} \longrightarrow 0$, almost uniformly, where $S_n(A):=\sum_{k=1}^n a_{n,k}x_k$.
        \item[(b)] If $\alpha=\tfrac12$, $\displaystyle \limsup_{n\to\infty} \frac{S_n(A)}{b_n} \underset{b.a.u.}{\leq} \sqrt2 K \ell$.
    \end{enumerate}
\end{mainthm}

Thus, our result extends the general framework of Li and Tomkins to the non-commutative setting. It is not limited to the special case $\sqrt{2n\log n}.$ The constant $\sqrt{2}$ is the key point of this result. Under the corresponding assumptions, our bound improves the constant $6\sqrt{2}$ in \cref{eqn:Li tomkins LL bound} to $\sqrt 2$. This improvement is also important in the classical setting. In the finite-variance logarithmic case, take $B(n)=\sqrt{2n\log n}.$ Then the truncated second-moment parameter satisfies $\ell^2=\frac{1}{2}\mathbb E(X^2).$ Hence, part $(b)$ of \Cref{mainthm:LIL deterministic weights} reduces to $\limsup_{n\to\infty} \frac{|S_n(A)|}{\sqrt{2n\log n}} \leq K\lrp{\mathbb E X^2}^{1/2} ~\text{a.s.}$ Thus, in the classical finite-variance case, our bound agrees with the sharp bound of Chen and Hao \cite{chen2011remark}. At the same time, our theorem keeps the more general normalization $b_n=B(n)$. In the critical case $\alpha=1/2$, no finite second moment is required. For $\alpha\in(1/2,1)$, the finite second moment assumption can also be replaced by a weaker condition on the truncated second moment; see \Cref{rem:L2 assumption}. In this sense, the result extends the sharp classical bound to the non-commutative setting and to a much larger class of distributions.

One main difficulty in proving this result is that the usual non-commutative LIL method cannot be used here. In the earlier works of Konwerska and Zeng, and also in our companion paper, the argument is built around a martingale. This makes it possible to control the supremum over a block through the column tail probability $\operatorname{Prob}_c$, and then use non-commutative Doob and Chebyshev inequalities together with a Borel-Cantelli argument \cite{konwerska2008law,zeng2015kolmogorov,panja2026non}. 

In the present weighted problem, however,  $S_n(A)=\sum_{k=1}^n a_{n,k}x_k $ does not form a martingale, since the coefficients change with $n$. Thus, this standard maximal-inequality argument for partial sums cannot be applied directly in this setting. We overcome this problem by working directly with the definitions of almost uniform and bilaterally almost uniform boundedness. We construct suitable projections for the different parts of the sum and control their complements through trace estimates. This gives the required a.u. and b.a.u. bounds without embedding the whole weighted sequence into a single martingale.

The assumptions are motivated by the classical law of the logarithm. They involve the tail of $x_1$ and its truncated second moment. These conditions help us control both the large values and the main fluctuations of the sum. The critical case $\alpha=1/2,$ is more delicate. Here we use a truncation of the random operators. We split the sum into a main part and tail parts. The main part is controlled by the non-commutative exponential inequality in \Cref{thm: exp ineq JFA}, as carried out in \Cref{lem:y_n part}. The large part is controlled by tail estimates and a non-commutative Borel-Cantelli argument (cf. \Cref{lem:v_n part}). 

A second important tool is dyadic blocking. We consider blocks $2^{k-1}<n\leq2^k.$ Regular variation gives a comparison between $b_n$ and $B(2^k)$ on each block. This allows us to replace the original problem by a sequence of block problems.

For the intermediate part of the truncated sum, we use an argument inspired by Cuzick's method for weighted sums \cite{cuzick1995strong}. This part is more delicate in the non-commutative setting. 
The intermediate terms are neither going to $0$ in a.u. nor uniformly bounded by constant, and the truncation level of these terms changes with $n$.
Moreover, the scalar centering and Borel-Cantelli argument used by Cuzick cannot be applied directly to non-commuting variables. We overcome this by working on dyadic blocks and writing the intermediate sum as a centered finite martingale part plus a deterministic mean part; see \Cref{lem:w_n part}. The centered part is controlled by a non-commutative exponential inequality and a non-commutative analogue of Cuzick's lemma, (cf. \Cref{lem:NC_Cuzick_2.2}).

A further difficulty appears since index $\alpha=1/2$, where the natural regularly varying integrand has index $-1$ and we are unable to use Karamata's theorem \cite{binghamregular1989}. We avoid this critical case by a suitable comparison that changes the relevant index from $-1$ to $-3/2$; see \Cref{lem:cuzick's methodology}. This gives an almost-uniform bound for the intermediate part. By choosing the truncation parameter suitably, this bound can be made arbitrarily small; see \Cref{lem:w_n part}.

For $\alpha\in(1/2,1)$, we assume that $x_1\in L_2(\CM,\tau)$. This assumption is not needed in the classical result of Li and Tomkins \cite{li2003law}. In the non-commutative setting, however, using the decomposition in \cref{eqn:S_n(A) decompose}, we can conclude bilaterally almost uniform convergence to zero without assuming the $L_2$ condition. On the other hand, the $L_2$ assumption enables us to establish almost uniform convergence.

The paper is organized as follows. In \Cref{sec: prelim}, we collect the basic tools used in the proofs. We recall almost uniform convergence and boundedness, Doob's maximal inequality, and the non-commutative Chebyshev and Borel-Cantelli lemmas. In \Cref{sec: Weighted LIL for martingale difference sequence}, we prove the weighted LIL for martingale difference sequences. We first treat monomial weights. We then pass to polynomial weights by a blocking argument. Finally, we extend the result to suitable functions of bounded variation and power series by approximation. The self-adjointness assumption is removed by the standard $2\times2$ matrix argument. In \Cref{Sec:non-commutative Law of Logarithm for Deterministic Weighted Sums}, we study general bounded deterministic triangular arrays. We use regular variation, truncation, exponential estimates, and dyadic blocking. The case $\alpha>1/2$ and the critical case $\alpha=1/2$ are treated separately.

\section{Preliminaries}\label{sec: prelim}

In this section, we briefly discuss all the objects that are necessary for this article. Throughout this paper, unless otherwise mentioned, we fix a finite von Neumann algebra $\CM \subseteq \CB(\CH)$ equipped with a faithful, normal tracial state $\tau$. We will denote the $*$-algebra of closed, densely defined operators on $\CH$, which are affiliated to $\CM$, by $\tilde{\CM}$. These are considered as the non-commutative analogue of the classical measurable functions on a probability space. By the spectral theorem of the unbounded operators, the tracial state $\tau$ can be seen to be extended to the positive part of $\title{\CM}$. Using this for every $1\leq p\leq \infty$, one can define the non-commutative $L_p$-spaces equipped with the norm $\norm{x}_p:= \tau(\abs{x}^p)^{1/p}$, where we set $L_\infty(\CM):= \CM$. For every $a \in \tilde{\CM}$, let $\abs{a}= \int_0^\infty \lambda de_\lambda$ be the spectral decomposition. Then as in the classical setup, one can define an analogue of the decreasing rearrangement, $\mu_{(\cdot)}:[0, \infty) \to [0, \infty)$ by 
\begin{equation*}
    \mu_t(a):= \inf \{ s\geq 0: \tau(e_{(s, \infty)}(\abs{a})) \leq t \}.
\end{equation*}
It is known that $\mu_\cdot$ is decreasing, continuous from right and $a \in L_p(\CM)$ if and only if $\mu_{(\cdot)}(a) \in L_p((0, \infty))$ for all $a \in \tilde{\CM}$. For detailed exhibition on this topic, we refer to \cite{nelson1974notes,terp1981lp,fack1986generalized,hiai2021lectures}. We will use $\CP(\CM)$ to denote the lattice of projections of the von Neumann algebra $\CM$.

\begin{defn}\label{defn: auleq}
     Let $(x_n)_{n \in \N}$ be a sequence in $\widetilde{\CM}$ For $K\geq 0$, we say that $(x_n)_n$ is almost uniformly (a.u) bounded by $K$ if, for every $\varepsilon>0$, there exists a projection $e\in N$ with $\tau(1-e)<\varepsilon$ such that, for every $\delta>0$, 
    $\limsup_{n\to\infty}\|x_ne\|_\infty\leq K+\delta.$

    We say that $(x_n)_n$ is bilaterally almost uniformly (b.a.u) bounded by $K$ if 
    $\limsup_{n\to\infty}\|ex_ne\|_\infty\leq K+\delta.$  
    We write 
    $\limsup_{n\to\infty} x_n\auleq K$ for a.u boundedness and 
    $\limsup_{n\to\infty} x_n\bauleq K$ for b.a.u boundedness.
\end{defn}

\begin{defn}\label{defn: auconv}
    Let $(x_n)_{n\in\mathbb{N}}$ be a sequence in $\widetilde{\CM}$ and
    $x\in\widetilde{\CM}$. If for any $\epsilon>0$ there exists a projection
    $e\in\CP(\CM)$ with $\tau(1-p)<\epsilon$ such that    
    \begin{enumerate}
        \item[(a)] $\|(x_n-x)e\|_\infty \longrightarrow 0$ as $n\to\infty$,
        we say that $x_n$ converges to $x$ almost uniformly.
        
        \item[(b)] $\|e(x_n-x)e\|_\infty \longrightarrow 0$ as $n\to\infty$,
        we say that $x_n$ converges to $x$ bilaterally almost uniformly.
    \end{enumerate}
\end{defn}

We also recall the $\ell_\infty$-column valued non-commutative $L_p$-space.
\begin{align*}
    L_p(\CM; \ell_\infty^c):= \{ x=(x_n) \subset L_p(\CM, \tau): \text{ there exists } a \in L_{2p}(\CM, \tau) \text{ and } (y_n) \subset \CM \text{ such that } x_n=y_na\}.
\end{align*}
When $p\geq 2$, $L_p(\CM; \ell_\infty^c)$ becomes a Banach space under the norm defined by $\norm{x}_{L_p(\CM; \ell_\infty^c)}:= \inf \{(\sup_{n \in \N} \norm{y_n}_\infty)a \}$, $x \in L_p(\CM; \ell_\infty^c)$, where the infimum is taken over all such factorization of the element $x$. Then we have the following asymmetric Doob's maximal inequality,
\begin{thm}\label{doob}\cite{junge2002doob}
    Let $4 \le p \le \infty$. For every $L_p(\mathcal{M})$-bounded martingale
    $x=(x_n)_{n \in \N}$:
    $$\norm{(x_n)_{n \in \N}}_{L_p(\ell_\infty^c)}\leq 2^{2/p}
    \sup_{n\in \N}\|x_n\|_p.$$    
 \end{thm}
Next, we recall the definition of the column-wise tail probability of a sequence $(x_n)$ in $\widetilde{\CM}$ and a few of its properties and the associated non-commutative analogues of Chebychev's inequality and Borel-Cantelli lemma. See \cite{konwerska2008law} and  \cite{zeng2015kolmogorov}. For any $t>0$, define
$$
\operatorname{Prob}_c\left( \sup_{n \in \N} \|x_n\| > t \right) 
:= \inf \left\{ s > 0 : \exists \text{ projection } e \text{ with } \tau(1 - e) < s \text{ and } \|x_n e\|_\infty \leq t \ \forall n \in \N \right\}.
$$
This satisfies:
$$
\operatorname{Prob}_c\left( \sup_{n \in \N} \|x_n\| > t \right) 
\leq 
\operatorname{Prob}_c\left( \sup_{n \in \N} \|x_n\| > r \right), \text{ whenever } t \geq r.
$$
 Moreover, if each coefficient $a_n \geq 1$ for $n \in \N$, then 
\begin{equation}\label{scaling}
\operatorname{Prob}_c\left( \sup_{n \in \N} \|x_n\| > t \right) 
\leq 
\operatorname{Prob}_c\left( \sup_{n \in \N} \|a_i x_n\| > t \right).
\end{equation}

\begin{lem}\cite[Lemma 3.16]{konwerska2008law}\label{lem: Chebyshev inequality lemma}
    Let $(x_n)_{n \in \N}$ be a sequence of random variables. For $t > 0$ and $1 \leq p < \infty$,
    $$
    \operatorname{Prob}_c\left( \sup_n \|x_n\| > t \right) \leq t^{-p} \|x\|_{L_p(\ell_\infty^c)}^p.
    $$
\end{lem}

\begin{lem}\cite[Lemma 3.17]{konwerska2008law}\label{lem: Borel-Cantelli lemma}
    Let $(k_n)_{n \geq n_0}$ be an increasing sequence of integers and $(z_n)$ be a sequence of random variables. Then 
    $$
    \forall \delta>0, ~~ \sum_{n \geq n_0} \operatorname{Prob}_c\left( \sup_{k_n\leq m< k_{n+1}} \|z_m\| > \gamma + \delta \right) < \infty
~~ \Longrightarrow~~ \limsup_{n \to \infty} z_n \underset{a.u.}{\leq } \gamma.
    $$
\end{lem}

We end this section stating the definition of independent random variables. Let $\CM_1, \CM_2$ be two subalgebras of $\CM$. 
\begin{defn}
    We say that $\CM_1$ and $\CM_2$ are independent if 
    $$
    \tau(xy)= \tau(x) \tau(y)
    $$
    for all $x \in \CM_1$ and $y \in \CM_2$. Moreover, a sequence of random variables $(x_i)_{i \geq 1}$ is said to be independent if for all $i\geq 1$ the von Neumann algebras  $\CM_i:= \{1, x_1,...,x_i \}''$ and ${x_{i+1}}''$ are independent.
\end{defn}

\section{Weighted LIL for martingale difference sequence}\label{sec: Weighted LIL for martingale difference sequence}

In this section, we prove \Cref{mainthm: polynomial bound tomkins}, which establishes a weighted law of the iterated logarithm for martingale difference sequences that generalizes \cite[Theorem 1.2]{panja2026non}. This theorem can be easily recovered from \Cref{lem: tomkins 1975 1st lemma} by choosing $B_n=s_n$ and $\alpha_n=1.$ As mentioned in the previous section, we fix a finite von Neumann algebra $\CM$ equipped with a faithful, normal tracial state $\tau$. We further assume an increasing filtration $(\CM_n)_{n \in \N}$ of von Neumann subalgebras of $\CM$ and for all $n \in \N$, $E_n: \CM \to \CM_n$ be the associated faithful, normal, $\tau$-preserving conditional expectation, satisfying $E_m \circ E_n= E_{min\{m,n\}}$ for all $m,n \in \N$. We recall that a sequence $(x_n)$ in $\widetilde{\CM}$ is called a martingale sequence with respect to the filtration $(\CM_n)$ if we have $E_n(x_{n+1})= x_n$ for all $n \in \N$. We further have the martingale differences defined by $dx_n:= x_n- x_{n-1}$ for all $n \in \N$ (with $x_{-1}=0$). Without loss of generality, we assume that $\CM_0 = \C$ and $E_0$ denotes the trivial conditional expectation from $\CM$ onto $\CM_0:= \C.1$. We begin by recalling the following theorem from \cite[Theorem 2.10]{panja2026non}.

\begin{thm}\label{thm: exp ineq JFA}
    Let $ (x_k)_{k=0}^n $ be a finite self-adjoint martingale sequence and $ d_k := dx_k $ represent the associated martingale differences.  Assume that the following conditions are true:
\begin{itemize}
    \item[(i)] $x_0 = 0 $, 
    \item[(ii)] $ \|d_k\| \leq M $ for all $1 \leq k \leq n$,
   \item[(iii)] $ \sum_{k=1}^n E_{k-1}(d_k^2) \leq D^2 1 $.
\end{itemize}
 Then, we have the following:
\begin{enumerate}
    \item For all $\lambda>0$, 
    $$
    \tau\left( e^{\lambda x_n} \right) \leq \exp \left[ \frac{e^{\lambda M} - \lambda M -1}{M^2} D^2 \right]
    $$
    \item For any $ \epsilon \in (0, 1] $ and $ \lambda \in \left[ 0, \frac{3 \epsilon}{M} \right] $,
    $$
    \tau\left( e^{\lambda x_n} \right) \leq \exp\left[ \frac{1}{2} (1 + \epsilon) \lambda^2 D^2 \right].
    $$
\end{enumerate}
\end{thm}
The proof of the following result is close to \cite[Theorem 3.1]{panja2026non}. However, to study the weighted version of LIL, we have to substitute $s_n$ and $u_n$ in \cite[Theorem 3.1]{panja2026non} with proper coefficients. We include the proof here for the sake of completeness.

\begin{lem}\label{lem: tomkins 1975 1st lemma}
    Let $(x_n)$ be a  non-commutative self-adjoint martingale sequence in $\mathcal{M}$ with associated difference sequence $(d_n).$ Let $\{\alpha_n\}$, $\{B_n\}$, and $\{c_n\}$ be
    positive real sequences satisfying:
    \begin{enumerate}
        \item $B_n \uparrow \infty$,
        \item $B_{n+1} \sim B_n$,
        \item $c_n^2 \log\log B_n^2 \to 0$,
        \item for some $N > 0$,
        $$\|d_m\|_\infty \le c_n\alpha_nB_n ~\text{ and } ~ \sum_{m=1}^n E_{m-1}(d_m^2)\le \alpha_n^2 B_n^2 \cdot \mathbf 1.$$
        for all $n > N$ and all $1\le m\le n$.
    \end{enumerate}
    Then
    $$\limsup_{n \to \infty}
      \frac{x_n}{\lrp{2 B_n^2 \log\log B_n^2}^{1/2}} \underset{a.u.}{\le } \limsup_{n \to \infty} \alpha_n.$$
\end{lem}

\begin{proof}
    Let $\alpha=\limsup \alpha_n$ and assume that $\alpha<\infty.$ Choose  arbitrarily $\delta>\delta'>\delta''>\delta'''>\alpha$ and $1<c<\delta/\delta'$. Now choose $\varepsilon_0\in (0,1]$ such that $\frac{1}{1+\varepsilon_0} \cdot \brac{\frac{\delta''}{\delta'''}}^2>1.$ Let $n_k=\min\{m: B_m\ge c^k\}.$ Now notice that, $B_{n_k}\sim c^k$.\\
    Also set $t_n^2=2\log\log B_n^2.$ Then $t_{n_k}\sim t_{n_k+1}.$ Now, $\exists$ $K_0>0$ such that, if $k\ge K_0$ then
    \begin{enumerate}\label{conditions}
        \item\label{item: 1975 1} $n_k > n_{k-1} \ge N$
        
        \item\label{item: 1975 2} $\delta B_{n_{k-1}} t_{n_{k-1}} > \delta' B_{n_k} t_{n_k}$
        
        \item\label{item: 1975 3} $\frac{1}{1+\varepsilon_0} \cdot \brac{\frac{\delta''}{\delta'''}}^2 t_{n_k}^2 \ge 4$
        
        \item\label{item: 1975 4} $\alpha_{n_k} < \delta'''$
        
        \item\label{item: 1975 5} $\delta'' t_{n_k} c_{n_k} \le 3 \delta''' \varepsilon_0$
    \end{enumerate}
    First note  that, since $B_nt_n\uparrow \infty$ and $n_{k-1}<n\le n_k$, we have $\frac{t_{n} B_{n}}{t_{n_{k-1}} B_{n_{k-1}}}\ge 1.$
    Now notice that,
    \begin{align*}
        \operatorname{Prob}_c\left(\sup_{{n_{k-1}<n\le n_k}}\left\|\frac{x_n}{t_{n} B_{n}}\right\|>\delta\right)
        &\le \operatorname{Prob}_c\left(\sup_{{n_{k-1}<n\le n_k}}\left\|\frac{t_{n} B_{n}}{t_{n_{k-1}} B_{n_{k-1}}}\cdot \frac{x_n}{t_{n} B_{n}}\right\|>\delta\right)
        \\
        &\le \operatorname{Prob}_c\left(\sup_{{n_{k-1}<n\le n_k}}\|x_n\|>\delta B_{n_{k-1}}t_{n_{k-1}}\right)
        \\
        &\le  \operatorname{Prob}_c\left(\sup_{{n_{k-1}<n\le n_k}}\|x_n\|>\delta'B_{n_k}t_{n_k}\right),~\text{ using } \eqref{item: 1975 2}
        \\
        &\le \operatorname{Prob}_c\left(\sup_{{n_{k-1}<n\le n_k}}\|x_n\|> \delta'' t_{n_k}B_{n_k}\right),~\text{ Since } \delta''<\delta'
        \\
        &\le \operatorname{Prob}_c\left(\sup_{{n_{k-1}<n\le n_k}}\left\|\frac{x_n}{\alpha_{n_k} B_{n_k}}\right\|> \delta'' t_{n_k}/\delta'''\right),~\text{ using } \eqref{item: 1975 4}.\tag{\theequation}\label{eqn: Probc}
    \end{align*}
    Set $\varepsilon_k:=\delta'' t_{n_k}/\delta'''$ and $\gamma_k := \alpha_{n_k} B_{n_k}.$ We now apply \Cref{thm: exp ineq JFA} directly to the martingale $\brac{x_n/\gamma_k}_{1\le n\le n_k}.$ Notice that the martingale differences of $x_n/\gamma_k$ are $d_m/\gamma_k$ for $1\le m\le n_k.$ Now note that,
    $\norm{\frac{d_m}{\gamma_k}}_\infty = \frac{\norm{d_m}_\infty}{\alpha_{n_k} B_{n_k}}\le c_{n_k}:= M_k.$
    Also notice that, $\sum_{m=1}^{n_k} E_{m-1}\brac{\frac{d_m^2}{\gamma_k^2}} \le \mathbf 1.$
    Now applying \Cref{thm: exp ineq JFA}, we have
    \begin{equation}\label{eqn: trace exp ineq}
        \tau\brac{\exp\brac{\frac{\lambda x_{n_k}}{\alpha_{n_k} B_{n_k}}}}\le \exp\brac{\frac{(1+\varepsilon_0)\lambda^2}{2}},~\text{ for } \lambda\in [0,3\varepsilon_0/c_{n_k}].
    \end{equation}
    Using \Cref{lem: Chebyshev inequality lemma}, in continuation of \eqref{eqn: Probc} we have for $p\ge 4$ and $\lambda>0$
    $$\operatorname{Prob}_c\left(\sup_{{n_{k-1}<n\le n_k}}\left\|\frac{\lambda x_n}{\alpha_{n_k} B_{n_k}}\right\|> \lambda \varepsilon_k\right)\le \brac{\lambda \varepsilon_k}^{-p} \norm{\brac{\frac{\lambda x_n}{\alpha_{n_k} B_{n_k}}}}^p_{L^p(\ell^\infty_c)}.$$
    An application of \Cref{doob} yields,
    $$ \norm{\brac{\frac{\lambda x_n}{\alpha_{n_k} B_{n_k}}}}_{L^p(\ell^\infty_c)} \le 2^{2/p} \norm{\frac{\lambda x_{n_k}}{\alpha_{n_k} B_{n_k}}}_p.$$
    By the functional calculus inequality $|u|^p \leq p^p e^{-p}(e^u + e^{-u}),$ applied to $\dfrac{x_{n_k}}{\alpha_{n_k} B_{n_k}}$ and using \eqref{eqn: trace exp ineq}:
    $$\norm{\frac{\lambda x_{n_k}}{\alpha_{n_k} B_{n_k}}}_p^p \le p^p e^{-p}\brac{\tau\brac{e^{\frac{\lambda x_{n_k}}{\alpha_{n_k} B_{n_k}}}}+ \tau\brac{e^{-\frac{\lambda x_{n_k}}{\alpha_{n_k} B_{n_k}}}}}\le 4p^pe^{-p} \exp\brac{\frac{(1+\varepsilon_0)\lambda^2}{2}}.$$
    Combining we have,
    \begin{align*}
        \operatorname{Prob}_c\left(\sup_{{n_{k-1}<n\le n_k}}\left\|\frac{\lambda x_n}{\alpha_{n_k} B_{n_k}}\right\|> \lambda \varepsilon_k\right) &\le 8\brac{\lambda \varepsilon_k}^{-p} p^pe^{-p} \exp\brac{\frac{(1+\varepsilon_0)\lambda^2}{2}} \\
        &\le 8\brac{\frac{p}{\lambda e \varepsilon_k}}^p\exp\brac{\frac{(1+\varepsilon_0)\lambda^2}{2}}.
    \end{align*}
    Put $p=\lambda \varepsilon_k$. Thus aforesaid equation becomes,
    $$\operatorname{Prob}_c\left(\sup_{{n_{k-1}<n\le n_k}}\left\|\frac{\lambda x_n}{\alpha_{n_k} B_{n_k}}\right\|> \lambda \varepsilon_k\right) \le 8 \exp\brac{\frac{(1+\varepsilon_0)\lambda^2}{2} - \lambda \varepsilon_k}.$$
    Now set $\lambda = \frac{\varepsilon_k}{1+\varepsilon_0}.$ Using \eqref{item: 1975 3} and \eqref{item: 1975 5} it follows that $p\ge 4$ and $\lambda\le 3\varepsilon_0/c_{n_k}$. Since $\varepsilon_k=\delta'' t_{n_k}/\delta''',$ we have
    \begin{align*}
        \operatorname{Prob}_c\left(\sup_{{n_{k-1}<n\le n_k}}\left\|\frac{\lambda x_n}{\alpha_{n_k} B_{n_k}}\right\|> \lambda \varepsilon_k\right) & 
        \le 8 \exp\brac{-\brac{\frac{1}{1+\varepsilon_0} \cdot \brac{\frac{\delta''}{\delta'''}}^2 \log\log B_{n_k}^2}}\\
        & \sim 8 \brac{(2\log c)\cdot k}^{-\brac{\frac{1}{1+\varepsilon_0} \cdot \brac{\frac{\delta''}{\delta'''}}^2}},~\text{since } B_{n_k}\sim c^k.
    \end{align*}
     Since $\frac{1}{1+\varepsilon_0} \cdot \brac{\frac{\delta''}{\delta'''}}^2>1$ by assumption, we have $\displaystyle\sum_k \operatorname{Prob}_c\left(\sup_{{n_{k-1}<n\le n_k}}\left\|\frac{x_n}{t_{n} B_{n}}\right\|>\delta\right)<+\infty.$ Since $\delta>\alpha=\limsup \alpha_n$ is already chosen, it follows by \Cref{lem: Borel-Cantelli lemma} 
    $$ \limsup_{n\to\infty}\frac{x_{n}}{(2B_{n}^{2}\log\log B_{n}^{2})^{1/2}}\underset{a.u.}{\le } \limsup_{n\to\infty}\alpha_{n}.$$
\end{proof}

We now use this result to obtain our first non-commutative weighted law of iterated logarithm for martingale difference sequences, weighted by monomials.

\begin{thm}\label{thm: matingale difference monomial with non-uniform variance}
    Let $(x_n)$ and $(d_n)$ are sequences as in \Cref{lem: tomkins 1975 1st lemma}. Further assume that $\tau(d_1)=0$, $E_{m-1}(d_m^2)=\sigma_m^2 \cdot \mathbf 1<\infty$ and $\max_{m \le n} \norm{d_m} \le M_n$ for all $m,n \in \N$, where $(M_n)\subset \R^+$. Fix $j\ge 0$ and define $B_n^2:= \sum_{m=1}^n m^{2j}\sigma_m^2$. Assume that for every $j\ge 0,$
    $B_n^2\to\infty~\text{ and } ~\frac{n^{2j} M_n^2 L(B_n^2)}{B_n^2}\to 0,~\text{ as } n\to \infty.$
    Then,
    \begin{equation}\label{eqn: martingale polynomial ineq}
        \limsup_{n\to\infty} \frac{\sum_{m=1}^n m^j d_m}{\sqrt{B_n^2 \log\log B_n^2}} \underset{a.u.}{\le } \sqrt 2
    \end{equation}
\end{thm}

\begin{proof}
    Fix $j\ge 0$ and define $x_n:=\sum_{m=1}^n m^j d_m.$ Note that $(x_n)$ is a sequence of self-adjoint martingales associated to the filtration $(\CM_n)$. Its martingale differences are $\Delta_m:=x_m-x_{m-1}=m^j d_m,~1\le m\le n.$
    Set $\alpha_n:=1$ and  $c_n:=\frac{n^jM_n}{B_n}.$ By hypothesis $\frac{n^{2j} M_n^2 L(B_n^2)}{B_n^2}\to 0$, therefore $c_n^2\log\log B_n^2\to 0$ as $n\to \infty.$\\
     Next, for every $1\le m\le n$, 
    $\|\Delta_m\|_\infty \le c_n B_n.$
    Since for all $m \in \N$, $E_{m-1}(d_m^2)=\sigma_m^2 \cdot \mathbf 1$, we obtain
    $\sum_{m=1}^nE_{m-1}(\Delta_m^2) = B_n^2\cdot \mathbf 1.$
    Also notice that, 
    \begin{align*}
        0\le B_{n+1}^2 -B_n^2 \le \|\Delta_{n+1}\|^2 \le (n+1)^{2j} M_{n+1}^2.
    \end{align*}
    Since by the assumption $\frac{n^{2j} M_n^2 L(B_n^2)}{B_n^2}\to 0$, we have $(B_{n+1}^2 -B_n^2)/B_{n+1}^2\longrightarrow 0$ as $n\to \infty.$ Hence, $B_{n+1}\sim B_n,~\text{ as } n\to \infty.$
    This follows that, all the conditions of \Cref{lem: tomkins 1975 1st lemma} are satisfied, by taking  $\alpha_n=1.$ Therefore, 
    $$ \limsup_{n\to\infty} \frac{x_n} {\sqrt{2B_n^2\log\log B_n^2}} \underset{a.u.}{\le} \limsup_{n\to\infty}\alpha_n = 1. $$

\end{proof}

\begin{cor}\label{cor: matingale difference monomial}
Let $(x_n)$ and $(d_n)$ are sequences as in \Cref{lem: tomkins 1975 1st lemma}. Assume that $\tau(d_1)=0$, $E_{m-1}(d_m^2)=1$ and $\max_{m \le n} \norm{d_m} \le M_n$, where $(M_n)$ is a sequence of positive reals such that $\lim\limits_{n \to \infty} L(n) M_n^2/n=0.$
    Then, for every $j\ge 0$,
    \begin{equation}\label{eqn: martingale polynomial ineq with uniform bound}
        \limsup_{n\to\infty} n^{-j}(2n\log\log n)^{-1/2} \sum_{m=1}^n m^j d_m \underset{a.u.}{\le } (2j+1)^{-1/2}
    \end{equation}
\end{cor}

\begin{proof}
    The result follows immediately from \Cref{thm: matingale difference monomial with non-uniform variance} by taking $\sigma_m=1$ for all $m\ge 1.$ In this case,
    $B_n^2= \sum_{m=1}^n m^{2j}.$
    Now use the facts 
    $B_n^2 \sim \frac{n^{2j+1}}{2j+1} ~\text{ and } ~ \log\log B_n^2\sim \log\log n $
    in \eqref{eqn: martingale polynomial ineq}, we obtain
    $$\limsup_{n\to\infty} n^{-j}(2n\log\log n)^{-1/2} \sum_{m=1}^n m^j d_m \underset{a.u.}{\le } (2j+1)^{-1/2}. $$
\end{proof}

Next, we focus on establishing the above results for more general weights, namely polynomials. We begin with the following technical lemma.

\begin{lem}\label{lem: Block sequence estimates}
    Let $c>1$ and $r_n=\sqrt{2n\log\log n},~~n\in \N$. Define a sequence $n_k=\lfloor c^k\rfloor+1$ for $k\ge 1$.
    Then there exists an integer $K\ge1$ such that for every $k\ge K$:
    $$\sup_{n_{k-1}<n\le n_k} \left| \frac1{r_n}-\frac1{r_{n_k}} \right| \le \frac{6(c-1)}{r_{n_k}}.$$
\end{lem}
\begin{proof}
    First, observe that $n_k\sim c^k.$
    Choose $K_0\ge 1$, such that $n_{k-1}<n_k$ for all $k\ge K_0$. In particular, for all $k\ge K_0$:
    \begin{equation}\label{eq:block-ratio}
        1 \leq \frac{n_k}{n_{k-1}} \leq c + 1 \leq 2c.
    \end{equation}
    Now using $\dfrac{\log\log n_k}{\log \log n_{k-1}}\to 1$ as $k\to \infty$, there exists $K_1>K_0$ such that for all $k\ge K_1$
    \begin{equation}\label{eq:loglog-ratio}
        \log\log n_k\le 2\log \log n_{k-1}.
    \end{equation}
    Now notice that, for every $k\geq K_1$, \eqref{eq:block-ratio} and \eqref{eq:loglog-ratio} imply $r_{n_k}^2\le 8cr_{n_{k-1}}^2.$ 
    Now fix $k\geq K_1$ and let $n_{k-1}<n\leq n_k$. Since $n \leq n_k$, we have
    \begin{align}
        \left|\frac1{r_n}-\frac1{r_{n_k}}\right|\leq\frac{r_{n_k}^2-r_{n_{k-1}}^2}{r_{n_{k-1}}^2r_{n_k}}.\label{eq:second-estimate}
    \end{align}
    To estimate the numerator, write  $r_{n_k}^2-r_{n_{k-1}}^2=2\brac{(n_k-n_{k-1})\log\log n_k + n_{k-1}\lrp{\log\log n_k-\log\log n_{k-1}}}.$ 
    Since, $n_k = \lfloor c^k\rfloor+1 \leq c^k+1, ~ n_{k-1}\geq c^{k-1},$ we obtain 
    $n_k-n_{k-1} \leq (c-1)n_{k-1}+1.$
    Now choose $K_2>K_1$ such that for all $k\ge K_2$ we have 
    \begin{equation}\label{eq:n-difference} 
        n_k-n_{k-1} \leq 2(c-1)n_{k-1}.
    \end{equation}
    Using \eqref{eq:loglog-ratio} and \eqref{eq:n-difference}, for all $k\ge K_2$ we obtain, 
    $(n_k-n_{k-1})\log\log n_k \leq 4(c-1) n_{k-1}\log\log n_{k-1}.$
    Now, by the mean value theorem applied to $g(x)=\log\log x$, there exists $\xi_k\in(n_{k-1},n_k)$ such that 
    $\log\log n_k-\log\log n_{k-1} = (n_k-n_{k-1})/ \xi_k\log\xi_k.$
    Since $\xi_k\geq n_{k-1}$, we have $ n_{k-1} \lrp{ \log\log n_k-\log\log n_{k-1} } \leq \frac{n_k-n_{k-1}} {\log n_{k-1}}\le 2(c-1) n_{k-1}\log\log n_{k-1},$
    for all $k\geq K$, where $K>K_2$ so large that $\frac{1}{\log n_{k-1}}\le \log\log n_{k-1}$ for all $k\ge K.$\\
    Combining the above estimates, we have
    $r_{n_k}^2-r_{n_{k-1}}^2 \le 6(c-1) r_{n_{k-1}}^2,~k\ge K,$
    and therefore, \eqref{eq:second-estimate} yields 
    $$ \left| \frac1{r_n} -\frac1{r_{n_k}} \right| \leq \frac{6(c-1)} {r_{n_k}}, ~~ n_{k-1}<n\leq n_k,~~ k\geq K. $$
\end{proof}

For the remainder of this section, let $(r_n)_n$ and $(n_k)_k$ be the sequences obtained in \Cref{lem: Block sequence estimates}.

\begin{lem}\label{lem: polynomial weight difference au}
     Let $p(t)=\sum_{j=0}^{l} a_j t^j $ be a polynomial on $[0,1]$ and $(d_n)$ be the sequence as in \Cref{thm: matingale difference monomial with non-uniform variance}. Assume that for each $j=0,1,\ldots,l$, 
    \begin{equation}\label{eqn: monomial bound in lemma assumption}
        \limsup_{n\to\infty} \frac{1}{n^j r_n} \sum_{m=1}^{n} m^j d_m \underset{a.u.}{\le} 1. 
    \end{equation}
    Then there exists a constant $C_p > 0$ such that 
    \begin{itemize}\label{eq: weight difference bound}
        \item[\rm{(i)}] $\limsup_{k\to \infty} \left(\frac{1}{r_n} - \frac{1}{r_{n_k}}\right) \sum_{m=1}^{n} p\left(\frac{m}{n}\right) d_m \underset{a.u.}{\le} 6(c-1)A,$
        \item[\rm{(ii)}] $\limsup_{k\to \infty}\frac{1}{r_{n_k}} \sum_{m=1}^{n} \lrp{ p(m/n) - p(m/n_k) } d_m
        \underset{a.u.}{\le} C_p(c-1),$
    \end{itemize}
    where $A=\sum_{j=0}^l \abs{a_j}.$
\end{lem}
\begin{proof}
    First notice that, $\sum_{m=1}^n p \left(\frac{m}{n}\right)d_m  = \sum_{j=0}^l \frac{a_j}{n^j}\sum_{m=1}^n m^j d_m.$ Then using \eqref{eqn: monomial bound in lemma assumption}, for each $j$, there exists $e_j \in \mathcal{P}(M)$ with $\tau(1-e_j) < \varepsilon/(l+1)$ and $K'_j \geq 1$ such that $\left\|\frac{1}{n^j r_n}\left(\sum_{m=1}^n m^j d_m\right) e_j\right\|_\infty \leq 2,$ for all $n \geq K'_j$.
    Set $e := \bigwedge_{j=0}^l e_j$. Then $\tau(1-e) \leq \varepsilon$ and for all $n \geq N := \max_j K'_j$:
    \begin{align}
        \left\|\left(\sum_{m=1}^n p \left(\frac{m}{n}\right)d_m\right) e
        \right\|_\infty 
        &\leq \sum_{j=0}^l |a_j| \cdot \frac{1}{n^j} 
        \left\|\left(\sum_{m=1}^n m^j d_m\right) e\right\|_\infty= 2r_n A. \label{eq:mono-bound}
    \end{align}
    Since $n_k\to \infty$, there exists $K>\max\{N,K_2\}$ such that $n_{K_2-1}\ge K.$ Therefore, using \Cref{lem: Block sequence estimates}, for all $n_{k-1} < n \leq n_k$ and $n \geq K$:
    \begin{align*}
       \norm{\brac{\left(\frac{1}{r_n} - \frac{1}{r_{n_k}}\right) \sum_{m=1}^{n} p\left(\frac{m}{n}\right) d_m}e}_\infty 
        \le 6(c-1)A,
    \end{align*}
    This concludes \rm{(i)}.\\
    Now observe that,
    \begin{align*}
        \frac{1}{r_{n_k}}\sum_{m=1}^n \left[p \left(\frac{m}{n}\right) 
        - p \left(\frac{m}{n_k}\right)\right] d_m 
        &= \frac{1}{r_{n_k}} \sum_{j=0}^l a_j \left(\frac{1}{n^j} 
        - \frac{1}{n_k^j}\right) \sum_{m=1}^n m^j d_m.
        \label{eq:C-expand}
    \end{align*}
    
    For $n_{k-1} < n \leq n_k$ and each $j \geq 1$, by the mean value theorem applied to $g(x)=x^{-j}$ and using \eqref{eq:n-difference}, we have $\left|\frac{1}{n^j} - \frac{1}{n_k^j}\right|\leq \frac{2j(c-1)}{n_{k-1}^j}$,
    for all $k \ge K$. For $j = 0$, we have $1/n^0 - 1/n_k^0 = 0$, so the $j=0$ term 
    vanishes. Using \eqref{eq:mono-bound} (with $n$ replaced by the same $n$, and $e$ as above), for each $j \geq 1$ and $n \leq n_k$, 
    $\frac{1}{n^j}\left\|\left(\sum_{m=1}^n m^j d_m\right)e\right\|_\infty  \leq 2r_{n_k}.$
    Hence for $j \geq 1$ and for all $k\ge K$:    
    $$
        \frac{1}{r_{n_k}}\abs{\frac{1}{n^j}-\frac{1}{n_k^j}}
        \norm{\left(\sum_{m=1}^n m^j d_m\right)e}_\infty 
        \leq \frac{1}{r_{n_k}} \cdot \frac{2j(c-1)}{n_{k-1}^j} \cdot 2n^j r_{n_k} 
        \leq \frac{4j(c-1)n_k^j}{n_{k-1}^j} 
        \leq 4j(c-1)(2c)^j,
    $$
    where we used $n_k/n_{k-1} \leq 2c$ from \eqref{eq:block-ratio}. Therefore, summing over $j$, we have
    \begin{equation}\label{eq:C-bound}
       \norm{\brac{\frac{1}{r_{n_k}}\sum_{m=1}^n \left[p \left(\frac{m}{n}\right) 
        - p \left(\frac{m}{n_k}\right)\right] d_m}e}_\infty 
        \leq \sum_{j=1}^l |a_j| \cdot 4j(c-1)(2c)^j 
        =: (c-1)C_p,
    \end{equation}
    for a finite constant $C_p=\sum_{j=1}^l 4j(2c)^j|a_j|$. This proves \rm{(ii)}.
\end{proof}

\begin{lem}\label{lem: block sup estimate}
    Let $p(t)=\sum_{j=0}^{l} a_j t^j $ be a polynomial on $[0,1]$ such that  $\norm{p}_2:= \lrp{\int_0^1 p(t)^2 dt}^{1/2}= 1$ and $(d_n)$ be the sequence as in \Cref{cor: matingale difference monomial}. Then
    \begin{equation}
        \limsup_{k\to \infty} \frac{1}{r_{n_k}} \sum_{m=1}^{n} p(m/n_k) d_m \underset{a.u.}{\le} 1.
    \end{equation}
\end{lem}
\begin{proof}
    Let us first denote $A_k(n) = \frac{1}{r_{n_k}} \sum_{m=1}^{n} p(m/n_k) d_m$, for $n_{k-1}<n\le n_k.$ For fixed $k$, the weights $p(m/n_k)$ are deterministic scalars, so $(A_k(n))_{1 < n \leq n_k}$ forms a self-adjoint martingale with respect to the filtration $(\mathcal{M}_n)$. We will show that, for any $\delta>0,~\sum_k\operatorname{Prob}_c \left(\sup_{n_{k-1} < n \leq n_k} \norm{A_k(n)} > 1 + \delta \right)<\infty.$

    \noindent
    Choose $\varepsilon_0\in(0,1)$ so that $(1+\delta)^2/(1+\varepsilon_0)>1.$
    Consider the martingale differences $\Delta_m=\dfrac{1}{r_{n_k}}p(m/n_k) d_m.$ Now notice that, $\|\Delta_m\|_\infty \le \lrp{\|p\|_\infty/r_{n_k}}\cdot M_{n_k} = M_k$, for all $1\le m\le n_k.$
    Also notice that, $\sum_{m=1}^{n_k} E_{m-1}(\Delta_m^2)=D_k^2.$ 
    Since $\lim_{k \to \infty} \frac{1}{n_k}\sum_{m=1}^{n_k} p(m/n_k)^2 =\|p\|_2,$ we have
    \begin{equation}\label{eqn: D_k^2 growth}
        D_k^2 \sim \frac{1}{2\log\log n_k},~\text{ as } k\to \infty.
    \end{equation}    
    Therefore, $(A_k(n))_{1 < n \leq n_k}$ satisfies all the conditions of \Cref{thm: exp ineq JFA}, hence for any  $\lambda\in [0,3\varepsilon_0/M_k],$
    \begin{equation}\label{eq:exp-bound}
        \tau\left(\exp(\lambda A_k(n_k))\right) \leq \exp\left(\frac{(1+\varepsilon_0)\lambda^2 D_k^2}{2}\right).
    \end{equation}
    Let us now estimate the column-wise tail probability. Let $q \geq 4$ (to be chosen). Take $\lambda>0$.
    \begin{align*}
        \operatorname{Prob}_c \left(\sup_{n_{k-1} < n \leq n_k} \norm{A_k(n)} > 1 + \delta \right) &=\operatorname{Prob}_c \left(\sup_{n_{k-1} < n \leq n_k} \norm{\lambda A_k(n)} > \lambda (1 + \delta) \right)\\
        &\leq (\lambda (1 + \delta))^{-q} \left\|(\lambda A_k(n))_{n_{k-1}<n\leq n_k}\right\|_{L^q(\ell^\infty_c)}^q\\
        &\le (\lambda (1 + \delta))^{-q}\brac{2^{2/q}}^q \norm{\lambda A_k(n_k)}_q^q.
    \end{align*}
    Since $A_k(n_k)$ is self-adjoint, applying functional calculus for the function $|u|^q \leq q^q e^{-q}(e^u + e^{-u}),~ u \geq 0$ we obtain
    \begin{equation}\label{eq:func-calc}
        \tau(\lambda |A_k(n_k)|^q) \leq q^q e^{-q} \left[\tau(e^{\lambda A_k(n_k)}) + \tau(e^{-\lambda A_k(n_k)})\right].
    \end{equation}
    Combining all we have,
    \begin{align}
        \operatorname{Prob}_c \left(\sup_{n_{k-1}<n\leq n_k}\norm{\lambda A_k(n)} > \lambda(1 + \delta)
        \right) \notag 
        &\leq (\lambda(1+\delta))^{-q} \cdot  \brac{2^{2/q}}^q \cdot q^q e^{-q} \cdot  
        2\exp \left(\frac{(1+\varepsilon_0)\lambda^2D_k^2}{2}
        \right) \notag \\
        &\leq 8 \left(\frac{q}{\lambda(1+\delta)e}\right)^q
        \exp \left(\frac{(1+\varepsilon_0)\lambda^2D_k^2}{2}\right).
        \label{eq:prob-bound}
    \end{align}
    Now we optimize by setting $q=\lambda(1+\delta),$ we get
    $$ \operatorname{Prob}_c \left(\sup_{n_{k-1}<n\leq n_k}\norm{A_k(n)} > 1 + \delta\right) \le 8 \exp \left(\frac{(1+\varepsilon_0)\lambda^2D_k^2}{2}-\lambda(1+\delta)\right).$$
    Since $D_k\to 0,$ choose $K_0\ge 1$ such that $\dfrac{(1+\delta)^2}{(1+\varepsilon_0)D_k^2}\ge 4$ for all $k\ge K_0.$ Also recall from the assumption that $\lim_{n \to \infty} \frac{L(n) M_n^2}{n}=0.$ Therefore, $\dfrac{M_k}{D_k^2}\to 0$ as $k\to \infty$. Now choose $K\ge K_0$ such that
    $\frac{M_k}{D_k^2}\le \frac{3\varepsilon_0 (1+\varepsilon_0)}{1+\delta},~\forall k\ge K.$
    Now set $\lambda = \dfrac{1+\delta}{(1+\varepsilon_0)D_k^2}.$ Then it follows that $q\ge 4$ and $0\le \lambda\le \frac{3\varepsilon_0}{M_k}$. Now putting the value of $\lambda$ we have
    \begin{align*}
        \operatorname{Prob}_c \left(\sup_{n_{k-1}<n\leq n_k}\norm{A_k(n)} > 1 + \delta\right) &\le 8\exp\brac{-\frac{(1+\delta)^2}{2(1+\varepsilon_0)D_k^2}}\\
        &\sim  8\exp\brac{-\frac{(1+\delta)^2}{1+\varepsilon_0}\log\log n_k},~\text{ using \eqref{eqn: D_k^2 growth}}\\
        &\sim 8 \brac{k \log c}^{-\frac{(1+\delta)^2}{1+\varepsilon_0}}.
    \end{align*}
    By the choice of $\varepsilon_0,$ it follows that $\displaystyle\sum_k\operatorname{Prob}_c \left(\sup_{n_{k-1}<n\leq n_k}\norm{A_k(n)} > 1 + \delta\right)<\infty.$
    Hence using \Cref{lem: Borel-Cantelli lemma} $\limsup_{n\to\infty} A_k(n) \underset{a.u.}{\le } 1.$
\end{proof}

\begin{proof}[Proof of \Cref{mainthm: polynomial bound tomkins}]
    For any function $g$ in $[0,1]$ we define
    \begin{equation}\label{eqn: Z_n function}
        Z_n(g):=\frac{1}{r_n}\sum_{m=1}^n g\left(\frac{m}{n}\right)d_m.
    \end{equation}
    Now if $\|p\|_2=0$, we must have $p\equiv 0.$ Then \eqref{eqn: polynomial eqn} is trivial. So, without loss of generality, we may assume that $\|p\|_2=1.$ Also assume that,  $p(t)=\sum_{j=0}^l a_jt^j$ and choose $\varepsilon>0.$ 
    Now, let us define $A= \sum_{j=0}^l|a_j|$. Let $c>1$. Define a sequence $n_k=[c^k]+1,~k\ge 1.$ Now write, 
    \begin{equation}\label{eq: Z_n decomp}
        Z_n(p) 
        = \underbrace{\frac{1}{r_{n_k}} \sum_{m=1}^{n} 
            p\left(\frac{m}{n_k}\right) d_m}_{=: A_k(n)}
        + \underbrace{\left(\frac{1}{r_n} - \frac{1}{r_{n_k}}\right) 
            \sum_{m=1}^{n} p\left(\frac{m}{n}\right) d_m}_{=: B_k(n)}
        + \underbrace{\frac{1}{r_{n_k}} \sum_{m=1}^{n} 
            \left[p\left(\frac{m}{n}\right) - p\left(\frac{m}{n_k}\right)\right] 
            d_m}_{=: C_k(n)}.
    \end{equation}
    Using \Cref{cor: matingale difference monomial}, for each $j \in \{0, 1, \ldots, l\}$ we have,
    \begin{equation}\label{eq: monomial bound}
        \limsup_{n \to \infty} \frac{1}{n^j r_n} \sum_{m=1}^{n} m^j d_m \underset{a.u.}{\le } (2j+1)^{-1/2} \leq 1.
    \end{equation}
    Now using \Cref{lem: polynomial weight difference au} it follows that,
    $\limsup_{k\to \infty} \brac{B_k(n)+C_k(n)}\underset{a.u.}{\le } \brac{C+6A}(c-1),$
    for some constant $C> 0$. Now choose $c > 1$ close enough to $1$ that $(c-1)(6A + C) < \varepsilon.$ Combining the above steps with \Cref{lem: block sup estimate} it follows that, $\limsup_{n\to\infty} \brac{A_k(n)+B_k(n)+C_k(n)} \underset{a.u.}{\le } 1+\varepsilon.$
    Since $\varepsilon>0$ chosen arbitrarily, it follows from \eqref{eq: Z_n decomp}
    $\limsup_{n\to\infty} Z_n(p) \underset{a.u.}{\le } 1.$
\end{proof}

\begin{thm}\label{thm: self-adjoint trick}
    Let $(x_n)$ be a  non-commutative martingale sequence in a finite von Neumann algebra $\mathcal{M}$ with associated difference sequence $(d_n)$ satisfying $\tau(d_1)=0$, $E_{n-1}(d_n^* d_n)=E_{n-1}(d_n d_n^*)=1$ and $\max_{m \le n} \norm{d_m} \le M_n$, where $(M_n)$ is a sequence of positive reals such that $\lim\limits_{n \to \infty} L(n) M_n^2/n=0.$ 
	Let $p$ be a polynomial on $[0,1]$, then
    \begin{align}\label{eqn: polynomial eqn general}
        \limsup_{n \to \infty} \frac{1}{(\sqrt{2n L(n)})} \sum_{m=1}^n p(m/n) d_m \underset{a.u.}{\le } \|p\|_2.
    \end{align}
\end{thm}
\begin{proof}
    We employ the standard $2\times 2$ matrix trick, as in \cite[Theorem 3.3]{panja2026non}, to reduce the non-self-adjoint case to the self-adjoint setting. Let
    $\widetilde{\mathcal M}=\mathcal M \overline{\otimes} M_2(\mathbb C), ~\widetilde{\mathcal M}_k=\mathcal M_k \overline{\otimes} M_2(\mathbb C),~ k\ge 1.$
    For each $n\ge 0$, define
    $$\widetilde d_n = \begin{pmatrix} 0 & d_n\\ d_n^* & 0 \end{pmatrix} \in \widetilde{\mathcal M}.$$ 
    Then $\widetilde d_n$ is self-adjoint, 
    Since $(d_n)$ is a martingale difference sequence with respect to $(\mathcal M_n,E_n)$, it follows immediately that $(\widetilde d_n)$ is a self-adjoint martingale difference sequence with respect to $(\widetilde{\mathcal M}_n,\widetilde E_n)$, where for all $n \ge 0$
    \begin{align*}
        \widetilde{E}_n(\widetilde{d}_{n+1})= \begin{pmatrix} 0 & E_n(d_{n+1})\\ E_n(d_{n+1}^*) & 0 \end{pmatrix}.
    \end{align*}
    Also, for $n\ge 1$, 
    $ \|\widetilde d_n\|_\infty = \|d_n\|_\infty \le M_n. $
    Furthermore, notice that
    $$\widetilde d_n^2 =\begin{pmatrix} d_nd_n^* & 0\\ 0 & d_n^*d_n \end{pmatrix}, \text{ so }~ \widetilde E_{n-1}(\widetilde d_n^2)= \begin{pmatrix} \mathbf 1 & 0\\ 0 & \mathbf 1 \end{pmatrix} = 1_{\widetilde{\mathcal M}}.$$
    Therefore, all assumptions of \Cref{mainthm: polynomial bound tomkins} hold for $(\widetilde d_n)$. Hence,
    $$\limsup_{n\to\infty} \frac{1}{\sqrt{2nL(n)}} \sum_{m=1}^{n} p \left(\frac{m}{n}\right)\widetilde d_m \underset{a.u.}{\le} \|p\|_2~~ \text{in $\widetilde{\mathcal M}$.}$$
    Observe that $$\sum_{m=1}^{n} p \left(\frac{m}{n}\right)\widetilde d_m = \begin{pmatrix} 0 & \displaystyle\sum_{m=1}^{n} p \left(\frac{m}{n}\right)d_m \\[2ex] \displaystyle\sum_{m=1}^{n} p \left(\frac{m}{n}\right)d_m^* & 0 \end{pmatrix}.$$ 
    Let $\widetilde r = \begin{pmatrix} 0 & 0\\ 1 & 0 \end{pmatrix} \in \widetilde{\mathcal M}.$ Then $$ \widetilde r \left( \sum_{m=1}^{n} p \left(\frac{m}{n}\right)\widetilde d_m \right) = \begin{pmatrix} 0 & 0\\ 0 & \displaystyle\sum_{m=1}^{n} p \left(\frac{m}{n}\right)d_m \end{pmatrix}. $$
    Consequently,
    $$\limsup_{n\to\infty} \frac{1}{\sqrt{2nL(n)}}\widetilde r \left( \sum_{m=1}^{n} p \left(\frac{m}{n}\right)\widetilde d_m \right) \underset{a.u.}{\le} \|p\|_2.$$
    Since almost uniform boundedness does not depend on the ambient von Neumann algebra (see \cite[Remark 2.3]{panja2026non}) we obtain,
    $$ \limsup_{n \to \infty} \frac{1}{(\sqrt{2n L(n)})} \sum_{m=1}^n p(m/n) d_m \underset{a.u.}{\le } \|p\|_2.$$
\end{proof}

\begin{lem}\label{lem: NC tomkins lemma 2}
    Let $(d_n)_{n \geq 1}$ be a sequence of self-adjoint random variables in $(\mathcal{M}, \tau)$. Let $\{b_n\}$ be a non-decreasing sequence of positive real numbers such that
    \begin{equation}\label{eqn: LIL hypothesis}
        \limsup_{n \to \infty} \frac{1}{b_n} \sum_{m=1}^n d_m\underset{a.u.}{\leq} 1.
    \end{equation}
    Let $\{a_{nm}\}_{m,n \geq 1}$ be a double sequence of real numbers such that
    $a_{nm}=0, ~ \text{ whenever } m>n.$
    Assume moreover that, there exists $L \in \mathbb{R}$ such that 
            $\lim_{n \to \infty} a_{nm} = L ~ \text{for every fixed } m \geq 1.$
    Then for any positive sequence $\{\beta_n\}$ with $\liminf_{n \to \infty} \beta_n > 0$,
    $$\limsup_{n \to \infty} \frac{1}{\beta_n} \sum_{m=1}^{n} a_{nm} d_m  \underset{a.u.}{\leq} \limsup_{n \to \infty} \frac{b_n}{\beta_n}\sum_{m=1}^{n} |a_{nm} - a_{n,m+1}|.$$
\end{lem}

\begin{proof}
    Define $x_n:= \sum_{m=1}^n d_m$ and let $\varepsilon, \delta > 0$ be arbitrary. By \eqref{eqn: LIL hypothesis}, there exists a projection $e \in \mathcal{M}$ with $\tau(1 - e) < \varepsilon$ and $N = N(e) \in\mathbb{N}$ such that
    \begin{equation}\label{eqn: proj bound}
    \|x_m e\| \leq (1 + \delta) b_m ~~ \text{for all } m \geq N.
    \end{equation}
    Since $a_{nm}$ are real scalars, we have $\sum_{m=1}^{n} a_{nm} d_m 
        = \sum_{m=1}^{n} (a_{nm} - a_{n,m+1}) x_m,$  for $n\ge 1$,
    where we used $x_0 = 0$ and $a_{n,n+1}=0$. Now notice that, since $b_n$ is non-decreasing
    \begin{align}
        \norm{\brac{\sum_{m=1}^{n} a_{nm} d_m}e}_\infty &\le \sum_{m=1}^n \abs{a_{nm} - a_{n,m+1}}\cdot \norm{x_me}_\infty\notag\\
        &\le R_N + (1+\delta)\sum_{m=N}^n \abs{a_{nm} - a_{n,m+1}} b_m,~\text{ using \eqref{eqn: proj bound}}\notag\\
        &\le R_N + (1+\delta)b_n \sum_{m=N}^n \abs{a_{nm} - a_{n,m+1}},\label{eqn: decompose eqn}
    \end{align}
    where $R_N=\sum_{m=1}^{N-1} \abs{a_{nm} - a_{n,m+1}}\cdot \norm{x_me}_\infty.$ Now observe that, since $\lim_{n \to \infty} a_{nm} = L ~~ \text{for every fixed } m \geq 1$ and $\liminf_{n \to \infty} \beta_n > 0$, we have $\lim_{n\to \infty} \abs{a_{nm} - a_{n,m+1}}/\beta_n=0,$ for every fixed $m \geq 1$ implies $\lim_{n\to \infty}R_N/\beta_n=0.$
    It follows from \eqref{eqn: decompose eqn}, 
    $$\limsup_{n\to \infty} \beta_n^{-1}\norm{\brac{\sum_{m=1}^{n} a_{nm} d_m}e}_\infty \le (1+\delta) \limsup_{n\to \infty} \beta_n^{-1}b_n \sum_{m=N}^n \abs{a_{nm} - a_{n,m+1}}.$$
    This completes the proof of the lemma.
\end{proof}

\begin{lem}\label{lem: bounded variation bound}
    Let $(\mathcal{M}, \tau)$ and $(d_n)$ be same as in \Cref{lem: NC tomkins lemma 2} and  $\limsup\limits_{n \to \infty} \frac{1}{(2n \log\log n)^{1/2}}\sum_{m=1}^n d_m\underset{a.u.}{\leq} 1.$
    Then for any function $g : [0,1] \to \mathbb{R}$ of bounded variation with total variation $V_g$,
    \begin{equation}\label{eqn: NC lemma 3 conclusion}
        \limsup_{n \to \infty} \frac{1}{(2n\log\log n)^{1/2}}\sum_{m=1}^{n} g(m/n) d_m \underset{a.u.}{\leq} V_g + |g(1)|.
    \end{equation}
\end{lem}

\begin{proof}
    Define the double sequence of real numbers $\{a_{nm}\}_{m,n \geq 1}$ by
    \begin{equation}\label{eqn: weights def}
        a_{nm} := \begin{cases} 
                        g(m/n) & \text{if } 1 \leq m \leq n, \\ 
                        0 & \text{if } m > n.
                  \end{cases}
    \end{equation}
    Set $b_n = \beta_n := (2n\log\log n)^{1/2}$. For each fixed $m \geq 1$, as $n \to \infty$, we have $a_{nm} = g(m/n) \longrightarrow g(0).$ So condition of \Cref{lem: NC tomkins lemma 2} holds with $L = g(0)$. Since $a_{nm} = 0$ for all $m > n$,  and $b_n = \beta_n = (2n\log\log n)^{1/2}$, the hypotheses of \Cref{lem: NC tomkins lemma 2} are satisfied. Therefore, we have
    \begin{equation}\label{eqn: after lemma 2}
        \limsup_{n\to\infty} \frac{1}{(2n\log\log n)^{1/2}}\sum_{m=1}^n g(m/n)d_m \underset{a.u.}{\leq} \limsup_{n\to\infty}\sum_{m=1}^n |a_{nm} - a_{n,m+1}|.
    \end{equation}
    Therefore, for every $n\ge 1$, $\sum_{m=1}^n |a_{nm} - a_{n,m+1}| \le V_g + |g(1)|.$
    Since this bound holds for every $n \geq 1$, then by \eqref{eqn: after lemma 2} we have
    $$ \limsup_{n\to\infty} \frac{1}{(2n\log\log n)^{1/2}}\sum_{m=1}^n g(m/n)d_m \underset{a.u.}{\leq} V_g + |g(1)|.$$
\end{proof}

\begin{cor}\label{cor: weight by fns of bdd var}
    Let $(\CM, \tau)$ and $(d_n)$ be same as \Cref{thm: self-adjoint trick}. Let $f$ be a function on $[0,1]$ such that a sequence $(p_m)$ of polynomials exists such that $p_m(1)\to f(1)$ and $V_{p_m-f}\to 0$ as $m\to \infty$. Then
    \begin{align}\label{eqn: main eqn cor}
        \limsup_{n \to \infty} \frac{1}{(\sqrt{2n L(n)})} \sum_{m=1}^n f(m/n) d_m  \underset{a.u.}{\le } \|f\|_2.
    \end{align}
    In particular, \eqref{eqn: main eqn cor} holds if $f(t)=\sum_{j=0}^\infty a_jt^j$ is a power series with either
    \begin{itemize}
        \item[\rm{(i)}] radius of convergence greater than $1$, or
        \item [\rm{(ii)}] $\limsup\limits_{n\to \infty}\abs{a_n}^{1/n} = 1,~\sum_{j=0}^\infty \abs{a_j}<\infty$.
    \end{itemize}
\end{cor}
\begin{proof}
    First notice that, for any $t\in [0,1]$ and $m\ge 1$, 
    $\abs{p_m(t)-f(t)}\le V_{f-p_m}+ \abs{f(1)-p_m(1)},~\text{ see \cite{apostol1957mathematical}}.$
    By the hypothesis, $V_{f-p_m}\xrightarrow{m\to \infty}0$ and $\abs{f(1)-p_m(1)}\xrightarrow{m\to \infty}0.$ Hence we have, $\abs{p_m(t)-f(t)} \to 0,~\text{ as } m\to \infty,$ for each $t\in [0,1]$.
    Therefore $p_m\to f$ uniformly, so $\|p_m\|_2\to \|f\|_2$ as $m\to \infty.$\\
    Fix $m \geq 1$. Using \eqref{eqn: Z_n function} we decompose,
    \begin{equation}\label{eqn: decomp}
        Z_n(f) = Z_n(f-p_m)+Z_n(p_m).
    \end{equation}
    Now recall from \Cref{thm: self-adjoint trick}, 
    \begin{equation}\label{eqn: Z_n(p_m) bound}
        \limsup_{n\to \infty} Z_n(p_m) \underset{a.u.}{\le} \|p_m\|_2.
    \end{equation}
     Since $f$ is a BV function and $p_m$ is a polynomial (hence also BV), applying \Cref{lem: bounded variation bound} to $f-p_m$ we have
     \begin{equation}\label{eqn: Z_n(f-p_m) bound}
         \limsup_{n\to \infty} Z_n(f-p_m)\underset{a.u.}{\le} V_{f-p_m} +\abs{f(1)-p_m(1)}.
     \end{equation}
     Now using \eqref{eqn: Z_n(p_m) bound} and \eqref{eqn: Z_n(f-p_m) bound} to \eqref{eqn: decomp} we get
     $\limsup\limits_{n\to \infty} Z_n(f) \underset{a.u.}{\le} V_{f-p_m} +\abs{f(1)-p_m(1)} + \|p_m\|_2.$
     Letting $m\uparrow \infty,$ \eqref{eqn: main eqn cor} obtains.
     
     Now, let $f(t)=\sum_{j=0}^\infty a_jt^j$ and suppose that either \rm{(i)} or \rm{(ii)} holds. For each $m\ge 1,$ define $p_m(t)= \sum_{j=0}^m a_jt^j$. Then $p_m(1)\to f(1)$ and $V_{p_m-f}\to 0$ as $m\to \infty$ \cite{tomkins1975iterated}. Consequently, \eqref{eqn: main eqn cor} holds.
\end{proof}

\begin{rem}
    \Cref{cor: weight by fns of bdd var} also holds if $(d_n)$ is a sequence of mean-zero i.i.d.\ self-adjoint random variables satisfying $\|d_n\|_2\le 1$ and
    $$
        \max_{m\le n}\|d_m\|_\infty\le M_n, \qquad \frac{L(n)M_n^2}{n}\longrightarrow0.
    $$
    Indeed, with respect to the natural filtration, $(d_n)$ is a martingale-difference sequence, and hence \Cref{cor: weight by fns of bdd var} applies.
\end{rem}

\section{non-commutative Law of Logarithm for Deterministic Weighted Sums}\label{Sec:non-commutative Law of Logarithm for Deterministic Weighted Sums}

In this section, we prove \Cref{mainthm:LIL deterministic weights}, which establishes a weighted analogue of the non-commutative law of iterated logarithm. We aim to investigate the almost sure asymptotic behavior of weighted sums of non-commutative random variables under some tail conditions on their distributions. Let $(x_k)_{k\ge 1}$ be a i.i.d sequence of mean zero self-adjoint random variables of $(\CM,\tau)$. Let $A = \{a_{n,k} : 1\le k \le n,\ n\ge 1\}$ be a triangular array of real numbers and set 
\begin{equation}\label{eq:classD}
    \CD := \Big\{ A=(a_{n,k}) : \; K(A):=\sup_{n,k}|a_{n,k}| < \infty \Big\}.
\end{equation}
Throughout, for $A\in\CD$ we write $K := K(A) = \sup_{n,k}|a_{n,k}| < \infty.$

To generalize the asymptotic estimate more than the law of iterated logarithm, we normalise the averages by a real sequence obtained from a regularly varying function. In the following we state the required definition and few of the properties which are useful in our case.

\begin{defn}
    A positive measurable function $B:(0,\infty)\to(0,\infty)$ is said to be regularly varying at infinity with index $\alpha\in\mathbb{R}$, written $B\in\mathrm{RV}_{\alpha},$ if, for every $\lambda>0$, $\lim_{t\to\infty} \frac{B(\lambda t)}{B(t)} =\lambda^\alpha.$
    Equivalently, $B$ can be represented in the form $B(t)=t^\alpha L(t),$ where $L$ is a slowly varying function at infinity, that is, 
    $\lim_{t\to\infty}\frac{L(\lambda t)}{L(t)}=1,~ \lambda>0.$
\end{defn}

Let $B:(0,\infty)\to(0,\infty)$ be positive, strictly increasing, and regularly varying at infinity with index $\alpha$, where $\tfrac12 \le \alpha < 1.$ Set $b_n := B(n), ~ n\ge 1.$ Since $B$ is regularly varying it admits a generalized inverse $B^{-1}$ on its range, that is, $B^{-1}(s):=\inf\{u>0:B(u)>s\},~~ s>0.$ We also use the following standard properties of regularly varying functions. For the proofs and more information, we refer to \cite{binghamregular1989}.
\begin{itemize}
    \myitem{(RV1)}{itm:RV1} if $\varphi\in\mathrm{RV}_\alpha$ then $\log\varphi(t)\sim\alpha\log t$ and if $\alpha>0$ then                        $\varphi(t)\to\infty$.
    \myitem{(RV2)}{itm:RV2} if $B$ is positive, strictly increasing, and $B\in\mathrm{RV}_\alpha$ with $\alpha\in(0,1)$, its generalized inverse $B^{-1}$ satisfies $B^{-1}\in\mathrm{RV}_{1/\alpha}$ and $B^{-1}(B(t))\sim t$.
    \myitem{(RV3)}{itm:RV3} if $\theta$ is slowly varying and $\lambda_t\to\lambda\in(0,\infty)$ as $t\to\infty$, then $\theta(\lambda_t t)/\theta(t)\to1$.
    \myitem{(RV4)}{itm:RV4} If $f\in\mathrm{RV}_\alpha$ and $g$ is eventually positive, measurable, $g(t)\to\infty$, and  $g\in\mathrm{RV}_\sigma$ with $\sigma>0$, then $f\circ g\in\mathrm{RV}_{\alpha\sigma}$.
    \myitem{(RV5)}{itm:RV5} [Karamata's tail theorem] If $h$ is positive, measurable, locally integrable on some $[t_0,\infty)$, and $h\in\mathrm{RV}_\alpha$ with $\alpha<-1$, then for every $\varepsilon>0$ there is $c_\varepsilon$ such that $$\int_c^\infty h(t)dt \le (1+\varepsilon)\frac{ch(c)}{-\alpha-1}~~\text{for all } c\ge c_\varepsilon.$$
    \myitem{(RV6)}{itm:RV6} If $f\in\mathrm{RV}_{\rho}$ be a positive regularly varying function with index $\rho\in\mathbb{R}$. Then, for every $\eta>0$ and $A>1$, there exists $x_0=x_0(\eta,A)>0$ such that \[ \frac{f(y)}{f(x)} \le A\max\left\{ \left(\frac{y}{x}\right)^{\rho+\eta}, \left(\frac{y}{x}\right)^{\rho-\eta} \right\}, ~~~x,y\ge x_0. \]
\end{itemize}

\noindent Fix a $j\ge1$ and for $t>0$ let us define:
\begin{align}
    G(t) &:= \tau\left(1_{[t,\infty)}(|x_j|)\right) = \mu([t,\infty)), \label{eq:Gdef}\\
    H(t) &:= \tau\left(x_j^2 1_{[0,t)}(|x_j|)\right) = \int_{[0,t)} s^2 d\mu(s), \label{eq:Hdef}
\end{align}
where $\mu$ is the spectral distribution measure of $|x_1|$ with respect to $\tau$, defined by $ \mu(A):=\tau\lrp{1_A(|x_1|)},$ for every Borel subset $A$ of $[0,\infty)$. By identical distribution, $G$ and $H$ do not depend on the choice of $j$. For $j\ge1$ and $c>0$ we write the spectral truncation $x_j^{(c)} := x_j1_{[0,c)}(|x_j|),$ a bounded self-adjoint element of $\Wstar{x_j}$ with $\normi{x_j^{(c)}} \le c$.

\begin{example}
    Take $B(t)=\sqrt{t}\log t,$ for $t\geq e.$ Then, observe that $B\in\mathrm{RV}_{1/2}.$ Let $(x_n)$ be mean-zero, identically distributed, self-adjoint such that
    $$ G(t)=
    \begin{cases}
        1, & 0\leq t\leq1,\\
        t^{-2}, & t>1.
    \end{cases}$$
    Then, $\tau(\abs{x_1})<\infty$ and $x_1\notin L_2(\CM,\tau)$. Also note that, $H(t)=2\log t,~ t\geq1.$ Hence, $\tau\!\left(B^{-1}(|x_1|)\right)<\infty,$ and $ \frac{B^{-1}(t\log t)}{\log t}G(t)\longrightarrow0.$ Further, $\frac{n\log n}{b_n^2} H\!\left(\frac{b_n}{\log n}\right)=1.$ Thus, $ \ell=1.$ Therefore, all the assumptions of \Cref{mainthm:LIL deterministic weights}(b) holds. In particular, the theorem applies even though $x_1\notin L_2(\CM,\tau).$ It gives $\limsup_{n\to\infty} \frac{S_n(A)}{\sqrt n\,\log n} \underset{b.a.u.}{\leq} \sqrt{2}\,K.$
\end{example}

\begin{lem}\label{lem:G H integral expressions}
    Let $G$ and $H$ be as in \cref{eq:Gdef,eq:Hdef}. Then the following identities hold: 
    \begin{enumerate} 
        \item[(i)] For every $c\ge0$, $$ H(c) = \int_{0}^{c}2tG(t)dt - c^{2}G(c). $$ 
        \item[(ii)] If $\int_{[c,\infty)} sd\mu(s)<\infty,$ then $$ \int_{[c,\infty)} sd\mu(s) = cG(c) + \int_{c}^{\infty}G(t)dt. $$
    \end{enumerate}
\end{lem}
\begin{proof}
    $(i)$ Using Tonelli's theorem, we have
    $$
    \int_0^c 2tG(t)dt=\int_0^c 2t\int_{[0,\infty)}1_{[0,s]}(t)d\mu(s)dt
    =\int_{[0,\infty)}\lrp{\int_0^c 2t\cdot 1_{[0,s]}(t)dt}d\mu(s).
    $$
    Now notice that, $\int_0^c 2t \cdot 1_{[0,s]}(t)dt=\min(s,c)^2$. Hence
    $$
    \int_0^c 2tG(t)dt= \int_{[0,c)} s^2d\mu(s) + \int_{[c,\infty)} c^2d\mu(s) = H(c)+c^2G(c).
    $$
    $(ii)$ Using Tonelli's theorem, we have
    $$
    \int_c^\infty G(t)dt=\int_c^\infty\int_{[0,\infty)}1_{[0,s]}(t)d\mu(s)dt
    =\int_{[0,\infty)}\lrp{\int_c^\infty 1_{[0,s]}(t)dt}d\mu(s).
    $$
    Since $\int_c^\infty 1_{[0,s]}(t)dt= \max\{s-c,0\}$, we have
    $$
    \int_c^\infty G(t)dt=\int_{[c,\infty)}(s-c)d\mu(s)=\int_{[c,\infty)}sd\mu(s)-cG(c).
    $$
\end{proof}

\begin{lem}\label{lem: property1 of RV func}
    Let $\varphi$ be regularly varying of index $\alpha\in\mathbb R$. If $(\lambda_t)_{t>0}$ is a family of positive reals with $\lambda_t\to\lambda\in(0,\infty)$ as $t\to\infty$, then 
    $ \frac{\varphi(\lambda_t t)}{\varphi(t)}\longrightarrow\lambda^\alpha. $
\end{lem}
\begin{proof}
    By definition of regular variation, write $\varphi(t)=t^\alpha\theta(t)$ with $\theta$ slowly varying.
    Then, we get
    $$
    \frac{\varphi(\lambda_t t)}{\varphi(t)}=\lambda_t^\alpha\cdot\frac{\theta(\lambda_t t)}{\theta(t)}.
    $$
    Since $x\mapsto x^\alpha$ is continuous at $\lambda>0$ and $\lambda_t\to\lambda$, we have $\lambda_t^\alpha\to\lambda^\alpha$. Since $\lambda_t\to\lambda\in(0,\infty)$, \ref{itm:RV3} gives $\theta(\lambda_t t)/\theta(t)\to1$, as $t\to \infty$.
\end{proof}

\begin{lem}\label{lem:b_n c_n relation}
   Let $b_n := B(n), ~ n\ge 1,$ where $B:(0,\infty)\to(0,\infty)$ be positive, strictly increasing, and regularly varying at infinity with index $\alpha \in [\tfrac12,1)$ satisfying condition \eqref{item:first condition} of \Cref{mainthm:LIL deterministic weights}. Fix $\delta>0$ and, for $n\ge2$, set $c_n:=\delta b_n/\log n,~n\ge 2$. Then $c_n\to\infty$ as $n\to\infty$, and
    \begin{enumerate}
        \item[(a)] $\log c_n\sim\alpha\log n$;
        \item[(b)] $c_n\log c_n\sim\alpha\delta b_n$;
        \item[(c)] $B^{-1}(c_n\log c_n)\sim(\alpha\delta)^{1/\alpha}n$.
    \end{enumerate}
\end{lem}
\begin{proof}
    The map $t\mapsto B(t)/\log t$ is regularly varying of index $\alpha>0$, hence by \ref{itm:RV1}, we have $\log(B(t)/\log t)\sim\alpha\log t\to\infty$, so $B(t)/\log t\to\infty$. Observe that, $b_n/\log n\to\infty$, and multiplying by the fixed constant $\delta$, we get $c_n=\delta b_n/\log n\to\infty$ as $n \to \infty.$

    \noindent$(a)$ By \ref{itm:RV1}, $\log B(n)\sim\alpha\log n$. Since $\log c_n=\log\delta+\log b_n-\log\log n$, and $\log\delta$ is a fixed constant while $\log\log n=o(\log n)$, substituting $\log b_n=\log B(n)\sim\alpha\log n$ gives $\log c_n\sim\alpha\log n$.

    \noindent$(b)$ Multiply $c_n=\delta b_n/\log n$ by $\log c_n$ and use (a), we obtain
    $$c_n\log c_n = \delta\frac{b_n}{\log n}\cdot\log c_n \sim \delta\frac{b_n}{\log n}\cdot\alpha\log n = \alpha\delta b_n.$$
    \noindent$(c)$ Set $\lambda_n:=c_n\log c_n/b_n$. Now by (b), $\lambda_n\to\lambda:=\alpha\delta\in(0,\infty)$. Since $B^{-1}\in\mathrm{RV}_{1/\alpha}$, we may write $ B^{-1}(s)=s^{1/\alpha}L_2(s), $ where $L_2$ is slowly varying. Now notice that, $\lambda_n\to\lambda\in(0,\infty)$ and \ref{itm:RV3} gives 
    $$ \frac{B^{-1}(\lambda_n b_n)}{B^{-1}(b_n)} = \lambda_n^{1/\alpha} \frac{L_2(\lambda_n b_n)}{L_2(b_n)} \longrightarrow \lambda^{1/\alpha} = (\alpha\delta)^{1/\alpha}. $$ 
    Now using \ref{itm:RV2}, we have $B^{-1}(b_n)=B^{-1}(B(n))\sim n$. Since $\lambda_nb_n=c_n\log c_n$ by definition of $\lambda_n$, combining above we  get $B^{-1}(c_n\log c_n)\sim(\alpha\delta)^{1/\alpha}n$.
\end{proof}

\begin{lem}\label{lem: h properties}
    Let $B$ be as in \Cref{lem:b_n c_n relation}. Define $h(t):=\log t/B^{-1}(t\log t)$ for $t>1$. Then
    \begin{enumerate}
        \item[(i)] $h$ is positive and measurable on some interval $[t_0,\infty)$;
        \item[(ii)] $h\in\mathrm{RV}_{-1/\alpha}$, and $-1/\alpha<-1$ since $\alpha\in(0,1)$;
        \item[(iii)] $h$ is locally integrable on $[t_0,\infty)$;
        \item[(iv)] with $c_n$, $\delta$ as in \Cref{lem:b_n c_n relation}, $\ c_nh(c_n)\sim K(\alpha,\delta)b_n/n$
        as $n\to\infty$, where $K(\alpha,\delta):=\alpha\delta(\alpha\delta)^{-1/\alpha}$.
    \end{enumerate}
\end{lem}
\begin{proof}
    $(ii)$ Notice that, the map $g(t):=t\log t$ is measurable. Observe that $g(t)\to\infty$, and belongs to $\mathrm{RV}_1$. By \ref{itm:RV4} applied with $f=B^{-1}\in\mathrm{RV}_{1/\alpha}$ and $g\in\mathrm{RV}_1$, notice that $t\mapsto B^{-1}(t\log t)$ is regularly varying with index $1/\alpha$. Since $\log t\in\mathrm{RV}_0$, the quotient $h$ is regularly varying with index $\rho_h:=0-1/\alpha=-1/\alpha$, and $\rho_h<-1$ because $\alpha<1$.
    
    \noindent$(i), (iii)$ Note that $B^{-1}$ is Borel measurable. Then, there exists $t_0>0$ such that $B^{-1}(t\log t)>0$ for all $t\ge t_0$. It follows that $h$ is positive and measurable on $[t_0,\infty)$, thus proving (i). Now suppose that $[t_1,t_2]\subset [t_0,\infty)$ be compact. Since both $B^{-1}$ and  the function $t\mapsto t\log t$ are non-decreasing, we have
    $$ 0<B^{-1}(t_1\log t_1) \le B^{-1}(t\log t) \le B^{-1}(t_2\log t_2) <\infty, \qquad t\in[t_1,t_2]. $$ Therefore, $ 0\le h(t)  \le \frac{\log t_2}{B^{-1}(t_1\log t_1)}, ~ t\in[t_1,t_2]. $ Thus $h$ is bounded on every compact subinterval of $[t_0,\infty)$, and hence $h$ is locally integrable on $[t_0,\infty)$, which proves $(iii)$.

    \noindent$(iv)$ By definition of $h$ and \Cref{lem:b_n c_n relation} (b),(c), we get 
    $$
    c_nh(c_n)=c_n\cdot\frac{\log c_n}{B^{-1}(c_n\log c_n)}\sim\frac{\alpha\delta b_n}{(\alpha\delta)^{1/\alpha}n}=\frac{b_n}{n}\cdot\alpha\delta(\alpha\delta)^{-1/\alpha}=K(\alpha,\delta)\frac{b_n}{n}.
    $$
\end{proof}

\begin{prop}\label{prop:integral s is zero}
    Let $b_n$ be same as in \Cref{lem:b_n c_n relation} and $\mu$ be as in \eqref{eq:Gdef}. Under Conditions~\eqref{item:first condition} and~\eqref{item:second condition} of \Cref{mainthm:LIL deterministic weights}, for every $\delta>0$,
    $$
    \frac{n}{b_n}\int_{[\delta b_n/\log n,\ \infty)} sd\mu(s) \longrightarrow 0 \qquad (n\to\infty).
    $$
\end{prop}

\begin{proof}
    Fix $\delta>0$ and set $c_n:=\frac{\delta b_n}{\log n}, ~n\ge 2.$ Then by \Cref{lem:b_n c_n relation}, there exists $N\ge1$ such that $c_n$ satisfies the hypothesis of Condition~\eqref{item:second condition} of \Cref{mainthm:LIL deterministic weights} for all $n\ge N$. By \Cref{lem:G H integral expressions}, the identity applies at $c=c_n$ gives
    $$ \frac{n}{b_n}\int_{[c_n,\infty)}sd\mu(s) = \underbrace{\frac{n}{b_n}c_nG(c_n)}_{=:T_1(n)} \;+\; \underbrace{\frac{n}{b_n}\int_{c_n}^\infty G(t)dt}_{=:T_2(n)}. $$ 
    We prove that  $T_1(n)\to0$ and $T_2(n)\to0$ as $n\to \infty.$
    By definition of $T_1$, 
    $$ T_1(n)=\delta\frac{n}{\log n}G(c_n) = \delta\cdot\frac{\alpha}{(\alpha\delta)^{1/\alpha}}\cdot\frac{B^{-1}(c_n\log c_n)}{\log c_n}G(c_n)\cdot(1+o(1)), $$ 
    using $n\sim(\alpha\delta)^{-1/\alpha}B^{-1}(c_n\log c_n)$ and $\log n\sim\log c_n/\alpha$, both from \Cref{lem:b_n c_n relation}(a),(c). Since $c_n\to\infty$, Condition~\eqref{item:second condition} of \Cref{mainthm:LIL deterministic weights} gives $B^{-1}(c_n\log c_n)G(c_n)/\log c_n\to0$. So, $T_1(n)\to0$ as $n \to \infty.$

    Let $h$ be as in \Cref{lem: h properties}. Now observe that, Condition~\eqref{item:second condition} of \Cref{mainthm:LIL deterministic weights} gives $G(t)/h(t)\to0$ as $t\to\infty$. Set $\eta(c):=\sup\limits_{t\ge c}\frac{G(t)}{h(t)}. $ Note that $\eta(c)\to0$ as $c\to\infty$, and $\int_c^\infty G(t)dt\le\eta(c)\int_c^\infty h(t)dt$ for all $c>1$. By \Cref{lem: h properties}(ii),(iii), $h$ satisfies the hypotheses of \ref{itm:RV5} with $\rho=\rho_h=-1/\alpha$, so for every $\varepsilon>0$ there is $c_\varepsilon$ such that 
    \begin{equation}\label{eqn: int h ineq}
        \int_c^\infty h(t)dt\le (1+\varepsilon)\frac{\alpha}{1-\alpha}ch(c), ~~c\ge c_\varepsilon,
    \end{equation}
    using $-\rho_h-1=1/\alpha-1=(1-\alpha)/\alpha$. Since $c_n\to\infty$, there is $N_\varepsilon$ with $c_n\ge c_\varepsilon$ for $n\ge N_\varepsilon$. Now combining with \eqref{eqn: int h ineq} at $c=c_n$, we have 
    $$ \int_{c_n}^\infty G(t)dt \le (1+\varepsilon)\frac{\alpha}{1-\alpha}\eta(c_n)c_nh(c_n),~~n\ge N_\varepsilon. $$
    By \Cref{lem: h properties}(iv), we get $c_nh(c_n)\sim K(\alpha,\delta)b_n/n$. Therefore,
    $$ T_2(n)=\frac{n}{b_n}\int_{c_n}^\infty G(t)dt \le (1+\varepsilon)\frac{\alpha}{1-\alpha}K(\alpha,\delta)\eta(c_n)(1+o(1)),~~n\ge N_\varepsilon. $$
    Using $\eta(c_n)\to0$, we conclude that $T_2(n)\to0$ as $n \to \infty.$
\end{proof}

\begin{lem}\label{lem:sum G finite}
    Assume that $\tau\lrp{B^{-1}(|x_1|)}<\infty$ and recall $G$ from \eqref{eq:Gdef}. Then, for every $\Delta>0$, $\sum_{k=1}^\infty G(\Delta b_k)<\infty$.
\end{lem}
\begin{proof}
    First take $\Delta=1.$ Let $s>0$ and define $N(s):=\#\left\{k\ge 1: B(k)\le s\right\}.$ Now notice that, $N(s)=\sum_{k\ge 1} 1_{[B(k),\infty)}(s).$ Since $B$ is increasing, we have $N(s)\le B^{-1}(s)+1.$ Using normality of $\tau$ and functional calculus on $N(s)$, we have
    $\sum_{k=1}^\infty G(b_k)\le \tau(B^{-1}(|x_1|))+\tau(1)<\infty,$
    by hypothesis. 
    
    Now take $\Delta>0$ in general. Define, $B_\Delta:= \Delta B.$ Observe that, $B_\Delta$ is positive, strictly increasing, and regularly varying with index $\alpha$, with $B_\Delta^{-1}(s)=B^{-1}(s/\Delta)$. Therefore, $\sum_kG(\Delta b_k)\le\tau(B^{-1}(|x_1|/\Delta))+\tau(1)$. It remains to prove that $\tau\lrp{B^{-1}\lrp{\frac{|x_1|}{\Delta}}}<\infty.$ Since $B^{-1}\in\mathrm{RV}_{1/\alpha},$ using \Cref{lem: property1 of RV func} we have, $\frac{B^{-1}(s/\Delta)}{B^{-1}(s)} \longrightarrow \Delta^{-1/\alpha}, \text{ as } s\to\infty.$ Hence, there exists $s_0\in \R$ and $C<\infty$ such that, $B^{-1}(s/\Delta) \le C B^{-1}(s), ~~\forall s\ge s_0.$ Also, notice that, for $0\le s<s_0$ we have $B^{-1}(s/\Delta)\le B^{-1}(s_0/\Delta).$ Therefore, we can write
    \begin{align*}
        \tau\lrp{B^{-1}\lrp{\frac{|x_1|}{\Delta}}} & = \tau\lrp{B^{-1}\lrp{\frac{|x_1|}{\Delta}} 1_{[0,s_0)}(|x_1|)} +\tau\lrp{B^{-1}\lrp{\frac{|x_1|}{\Delta}} 1_{[s_0,\infty)}(|x_1|)}\\
        & \le B^{-1}(s_0/\Delta) \tau\lrp{1_{[0,s_0)}(|x_1|)} + C \tau\lrp{B^{-1}(|x_1|)1_{[s_0,\infty)}(|x_1|) } < \infty.
    \end{align*}
    This completes the proof.
\end{proof}

For $k\ge 0$ set 
$\CM_k := \Wstar{1,x_1,\dots,x_k} \subseteq \CM, ~ \CM_0 := \mathbb{C}\cdot 1.$ 
Since $\tau$ is finite, normal, and faithful, each subalgebra $\CM_k$ admits a (unique, trace-preserving) normal conditional expectation $ E_k : \CM\longrightarrow \CM_k. $

\begin{lem}\label{lem:y_k d_{n,k} introduce}
    Let $(x_k)$ and $(a_{n,k})$ are same as in \Cref{mainthm:LIL deterministic weights}. Fix $n\ge1$ and $c>0$, and set, for $1\le k\le n$, 
    $$ x_k^{(c)} := x_k 1_{[0,c)}\left(|x_k|\right), ~~ d_{n,k} := a_{n,k}\lrp{x_k^{(c)} - \tau\lrp{x_k^{(c)}}}, ~~ y_k := \sum_{j=1}^k d_{n,j},~~ y_0:=0. $$ 
    Then
    \begin{enumerate} 
        \item[(a)] $E_{k-1}(d_{n,k}) = 0$ for every $1\le k\le n$; 
        \item[(b)] $(y_k)_{k=0}^n$ is a self-adjoint martingale with respect to the filtration $(\CM_k,E_k)_{k=0}^n$. 
    \end{enumerate}
\end{lem}
\begin{proof}
    $(a)$ Write $f(t) := a_{n,k}\lrp{t\mathbf 1_{(-c,c)}(t) - \tau(x_k^{(c)})}$. Note that $d_{n,k}=f(x_k)\in\Wstar{x_k}$. Since, $\CM_{k-1}=\Wstar{1,x_1,\dots,x_{k-1}}$ is independent of $\Wstar{x_k}$, we have 
    $\tau(z x_k^{(c)}) = \tau(z)\tau(x_k^{(c)}), ~ \forall z\in\CM_{k-1}.$
    The scalar operator $\tau(x_k^{(c)})\cdot 1\in \CM_{k-1}$ and satisfies 
    $\tau\lrp{z\cdot \tau(x_k^{(c)})\cdot 1} = \tau(zx_k^{(c)}),~~\forall z\in \CM_{k-1}. $
    By uniqueness of the conditional expectation and  faithfulness of $\tau$, we have $E_{k-1}(x_k^{(c)}) = \tau(x_k^{(c)})\cdot 1$. It follows that, 
    $E_{k-1}(d_{n,k})  = 0. $
    
    $(b)$ Notice that, $y_k$ are self-adjoint for every $k.$ For each $j\le k$, $d_{n,j}\in W^*(x_j)$. Now $d_{n,j}\in\CM_k$ for every $j\le k$. Now $y_k\in \CM_k$. Write $y_k=y_{k-1}+d_{n,k}$. By linearity of $E_{k-1}$, we have
    $E_{k-1}(y_k)=E_{k-1}(y_{k-1})+E_{k-1}(d_{n,k}).$
    Since $y_{k-1}\in\CF_{k-1}$, we have $E_{k-1}(y_{k-1})=y_{k-1}$. By part (a), $E_{k-1}(d_{n,k})=0$. So 
    $E_{k-1}(y_k)=y_{k-1}+0=y_{k-1}.$
\end{proof}

\begin{lem}\label{lem:condition verificartion}
    Let $y_k$ and $d_{n,k}$ be as in \Cref{lem:y_k d_{n,k} introduce}. For a positive sequence $(\delta_n)$, set $c=\delta_n b_n/\log n$, with $(\delta_n)$ and $A\in\CD$ with $K=\sup_{n,k}|a_{n,k}|$, we obtain
    \begin{enumerate}
        \item[(a)] $\normi{d_{n,k}} \le M_n := 2Kc$ for every $1\le k\le n$;
        \item[(b)] $E_{k-1}(d_{n,k}^2) = \tau(d_{n,k}^2)\cdot 1$ for every $1\le k\le n$;
        \item[(c)] $\displaystyle\sum_{k=1}^n E_{k-1}(d_{n,k}^2) \le D_n^2\cdot 1$, where $D_n^2 := nK^2 H\left(\frac{\delta_n b_n}{\log n}\right)$.
    \end{enumerate}
\end{lem}

\begin{proof}
    $(a)$ By construction $\normi{x_k^{(c)}}\le c = \delta_n b_n/\log n$. So, $\normi{x_k^{(c)}-\tau(x_k^{(c)})} \le 2c$, and $\normi{d_{n,k}} = |a_{n,k}|\cdot\normi{x_k^{(c)}-\tau(x_k^{(c)})} \le K\cdot 2c = M_n$.

    $(b)$ Let $g(t) := f(t)^2$, with $f$ as in the proof of \Cref{lem:y_k d_{n,k} introduce}, a bounded Borel function, so $d_{n,k}^2 = g(x_k) \in \Wstar{x_k}$. Thus, uniqueness of the conditional expectation forces $E_{k-1}(d_{n,k}^2) = \tau(d_{n,k}^2)\cdot 1$.

    $(c)$ By definition, $\tau(d_{n,k}^2) = a_{n,k}^2\lrp{\tau\lrp{(x_k^{(c)})^2} - \tau(x_k^{(c)})^2} \le a_{n,k}^2 \tau\lrp{(x_k^{(c)})^2} = a_{n,k}^2 H(c)$, using~\eqref{eq:Hdef} and identical distribution of the $x_k$. Summing over $k=1,\dots,n$ and using $A\in\CD$ (so $a_{n,k}^2\le K^2$), we have $$ \sum_{k=1}^n E_{k-1}(d_{n,k}^2) = \sum_{k=1}^n \tau(d_{n,k}^2)\cdot 1 \le \lrp{\sum_{k=1}^n a_{n,k}^2}H(c)\cdot 1 \le n K^2 H(c)\cdot 1 = D_n^2\cdot 1.$$
\end{proof}

\subsection{Proof of \Cref{mainthm:LIL deterministic weights}(a)}\label{subsec:part(a)} Fix $\Delta>0$ and $n\ge1$ and let $1\le k\le n$. Using same notation as \Cref{lem:y_k d_{n,k} introduce}, define 
$x_k^{(c_n)}:=x_k1_{[0,c_n)}(|x_k|),~z_{k}^{(c_n)}:=x_k1_{[c_n,\infty)}(|x_k|), $ 
so that $x_k=x_k^{(c_n)}+z_{k}^{(c_n)}$. Here, we take $c_n:=\Delta b_n.$ Define, for $0\le l\le n,$
\begin{equation}\label{eqn:y_l, m_l, v_l introduce}
    y_l:=\sum_{k=1}^l a_{n,k}\lrp{x_k^{(c_n)}-\tau(x_k^{(c_n)})},\qquad
    m_l:=\sum_{k=1}^l a_{n,k}\tau(x_k^{(c_n)}),\qquad
    v_l:=\sum_{k=1}^l a_{n,k}z_{k}^{(c_n)}.
\end{equation}
Then, exactly,
\begin{equation}\label{eqn:fixed fraction decomp}
    S_n(A)\ =\ y_n+m_n+v_n.
\end{equation}
By \Cref{lem:y_k d_{n,k} introduce} applied with $c=\Delta b_n>0$, the difference $d_{n,k}:=a_{n,k}\lrp{x_k^{(c_n)}-\tau(x_k^{(c_n)})}$ satisfies $E_{k-1}(d_{n,k})=0$ for $1\le k\le n$, and $\lrp{y_l:=\sum_{k=1}^l d_{n,k}}_{l=0}^n$ is a self-adjoint martingale. To prove \Cref{mainthm:LIL deterministic weights}(a), we will set $\sigma^2:=\tau(x_1^2)<\infty.$

\begin{lem}\label{lem:centering fixed fraction}
    For every fixed $\Delta >0$, $m_n/b_n \longrightarrow 0$ as $n \to \infty.$
\end{lem}

\begin{proof}
    Since $\tau(x_k)=0$ for every $k$, $\tau(x_k^{(c_n)})=-\tau(z_{k}^{(c_n)})$. Using $|\tau(y)|\le\tau(|y|)$ and identical distribution, 
    $|m_n|/b_n \le K \cdot \frac{n}{b_n}\tau\lrp{|x_1|1_{[\Delta b_n,\infty)}(|x_1|)}=K\frac{n}{b_n} \int_{[\Delta b_n,\infty)} s\, d\mu(s).$ By \Cref{prop:integral s is zero}, it follows that $\lim\limits_{n\to \infty}\frac{m_n}{b_n}=0.$ 
\end{proof}

\begin{lem}\label{lem:tail kill fixed fraction}
    For every fixed $\Delta >0$, $\dfrac{v_n}{b_n}\ \xrightarrow{\ a.u.\ }\ 0$ as $n\to \infty.$
\end{lem}
\begin{proof}
    Set $q_k:=1_{[\Delta b_k,\infty)}(|x_k|)$. Now by \Cref{lem:sum G finite}, we have $\sum_k\tau(q_k)<\infty$. Choose $N$ with $\sum_{k\ge N}\tau(q_k)<\varepsilon_1/2$. Let $p_v^{(1)}:=\bigwedge_{k\ge N}(1-q_k)$. This gives $\tau(1-p_v^{(1)})<\varepsilon_1/2$. Since $B$ is increasing, we have $1_{[\Delta b_n,\infty)}(|x_k|)\le q_k$, for $n\ge k\ge N$. It follows that, $x_k 1_{[\Delta b_n,\infty)}(|x_k|)p_v^{(1)}=0$.
    
    Now choose $R<\infty$, such that $\sum_{k<N}\tau(1_{[R,\infty)}(|x_k|))<\varepsilon_1/2$. Set $p_v^{(2)}:=\bigwedge_{k<N}1_{[0,R)}(|x_k|)$. Since $b_n\to \infty,$ there exists $n_0\in \N$, such that $\Delta b_n\ge R,~\forall n\ge n_0.$ Then for all $n\ge n_0,$ we have $x_k1_{[\Delta b_n,\infty)}(|x_k|)p_v^{(2)}=0.$ Now take, $p_v= p_v^{(1)}\wedge p_v^{(2)}.$ Then, $\tau(1-p_v)<\varepsilon_1$ and $v_np_v=0$ for all $n\ge n_0.$ In particular
    $\left\|\frac{v_n}{b_n}p_v\right\|_\infty=0~~\text{for all } n\ge n_0.$    This completes the proof.
\end{proof}

\begin{lem}\label{lem:exp tail fixed fraction}
    Fix $\theta\in(0,2\alpha-1)$. For every $\varepsilon>0$ there exists $\Delta>0$ such that 
    \begin{itemize}
        \item[(1)] $ \sum\limits_{n=1}^\infty \tau\lrp{1_{(\varepsilon b_n,\infty)}(|y_{n}|)} < \infty. $
        \item[(2)] $ \limsup\limits_{n\to\infty}\dfrac{y_n}{b_n}\ \underset{a.u.}{\le}\ \varepsilon. $
    \end{itemize}
\end{lem}
\begin{proof}
    \emph{(1):}
    Fix $\Delta>0$ and $n\in \N.$ Now recall from \eqref{eqn:y_l, m_l, v_l introduce}, $\lrp{y_l}_{l=0}^n$ is a martingale.  By \Cref{thm: exp ineq JFA} with $M_n=2K\Delta b_n,~D_n^2=K^2\sigma^2 n$ as in \Cref{lem:condition verificartion}, for every $\lambda>0$, we have 
    $\tau\lrp{e^{\lambda y_n}} \le \exp\lrp{\frac{e^{\lambda M_n}-\lambda M_n-1}{M_n^2}D_n^2}.$ 
    Set $t_n:=\varepsilon b_n$ and 
    \begin{equation}\label{eqn:u_n defn}
        u_n := \frac{t_nM_n}{D_n^2} = \frac{\varepsilon b_n\cdot 2K\Delta b_n} {K^2\sigma^2 n} = \frac{2\varepsilon\Delta}{K\sigma^2} \frac{b_n^2}{n}.
    \end{equation} 
    Choose $\lambda_n:=\log(1+u_n)/M_n$, so $e^{\lambda_nM_n}=1+u_n$. By functional calculus, we have
    $$ \tau\lrp{1_{(t_n,\infty)}(y_n)}\ \le \exp\lrp{-\lambda_nt_n+\frac{D_n^2}{M_n^2}\lrp{u_n-\log(1+u_n)}}. $$ 
    Now write $a_n:=t_n/M_n=\varepsilon/(2K\Delta )$ and $s_n:=D_n^2/M_n^2$. Note that $s_nu_n=a_n$. Then notice that, $-\lambda_nt_n+\frac{D_n^2}{M_n^2}\lrp{u_n-\log(1+u_n)}\le -a_n\lrp{\log(1+u_n)-1}.$ Hence 
    \begin{equation}\label{eqn:exp tail intermediate} 
        \tau\lrp{1_{(\varepsilon b_n,\infty)}(y_n)} \le \exp\lrp{-\frac{\varepsilon}{2K\Delta }\lrp{\log(1+u_n)-1}}. 
    \end{equation}
    Since $b_n^2/n\in\mathrm{RV}_{2\alpha-1}$ and $2\alpha-1>\theta>0$, we have $(b_n^2/n)/n^\theta\in\mathrm{RV}_{2\alpha-1-\theta}$. Hence, by \ref{itm:RV1} there are $n_1$ and $c_\theta>0$ with $b_n^2/n\ge c_\theta n^\theta$ for $n\ge n_1$. Now, using \cref{eqn:u_n defn}, we have 
    $u_n\ \ge \frac{2\varepsilon \Delta c_\theta}{K\sigma^2}n^\theta =:  C_{\varepsilon,\Delta,\theta} n^\theta, ~ n\ge n_1,$
     so $\log(1+u_n)\ge\theta\log n+\log C_{\varepsilon,\Delta,\theta}$ for $n\ge n_1$. Substituting into \cref{eqn:exp tail intermediate}, for $n\ge n_1,$
     \begin{align*}
        \tau\lrp{1_{(\varepsilon b_n,\infty)}(y_n)} \le & \exp \lrp{-\frac{\varepsilon}{2K\Delta}\lrp{\theta \log n + \log  C_{\varepsilon,\Delta,\theta}}}\\
        \le& C'_{\varepsilon,\Delta,\theta}\;n^{-\varepsilon\theta/(2K\Delta)}
    \end{align*}
    for a constant $C'_{\varepsilon,\Delta,\theta}= \exp\lrp{-\frac{\varepsilon}{2K\Delta}\lrp{\log C_{\varepsilon,\Delta,\theta}-1}}$ not depending on $n$.  Now choose $ \Delta := \frac{\varepsilon\theta}{4K} > 0, $ so that $\varepsilon\theta/(2K\Delta)=2>1$. Then $\sum_n\tau\lrp{1_{(\varepsilon b_n,\infty)}(y_n)}<\infty$. The same argument applied to $-y_n$ (with the same $M_n,D_n^2$) gives $\sum_n\tau\lrp{1_{(\varepsilon b_n,\infty)}(-y_n)}<\infty$.

    \emph{(2):}     For each $n$, let $e_n:=1_{[0,\varepsilon b_n]}(|y_n|)\in\CP(\CM)$. Then
    $\normi{y_ne_n}\le\varepsilon b_n$, and
    $\tau(1-e_n)=\tau\lrp{1_{(\varepsilon b_n,\infty)}(|y_n|)}$. By
    \emph{(1)}, $\sum_n\tau(1-e_n)<\infty$. Given $\delta>0$, choose
    $N\in \N$ with $\sum_{n\ge N}\tau(1-e_n)<\delta$. Now set
    $p_N:=\bigwedge_{n\ge N}e_n\in\CP(\CM)$. Then,
    $\tau(1-p_N)\le\sum_{n\ge N}\tau(1-e_n)<\delta$. Since $p_N\le e_n$, $\normi{y_np_N}\le\normi{y_ne_n}\le\varepsilon b_n$, for every $n\ge N$.
    So $\limsup_n\normi{(y_n/b_n)p_N}\le\varepsilon$. As $\delta>0$ was arbitrary, this
    proves the claim.
\end{proof}

\begin{rem}\label{rem:L2 assumption}
    The assumption $x_1\in L_2(\CM,\tau)$ in \Cref{mainthm:LIL deterministic weights}(a) is stronger than what is actually needed. It may be replaced by the following condition: there exists $\theta_0>0$ such that, for every $\Delta>0$, $$ \limsup_{n\to\infty}\frac{n^{1+\theta_0}}{b_n^2}H(\Delta b_n)<\infty. $$ To see this, recall that the finite second moment assumption enters the proof of \Cref{lem:exp tail fixed fraction} only through the estimate $ D_n^2\le K^2\sigma^2 n. $ If the condition above is assumed instead, one can use $ D_n^2\le K^2 n H(\Delta b_n), $ and the rest of the argument goes through in the same way. Consequently, the conclusion of \Cref{mainthm:LIL deterministic weights}(a) still holds under this condition, even when $x_1$ does not have a finite second moment.
\end{rem}

Now we will proof the part (a) of \Cref{mainthm:LIL deterministic weights}.

\begin{proof}[Proof of \Cref{mainthm:LIL deterministic weights}(a)]
   Fix $\theta\in(0,2\alpha-1)$. For $j\ge1$ set $\varepsilon_j:=2^{-j}$ and $\Delta _j:=\Delta $ as constructed in \Cref{lem:exp tail fixed fraction}. Fix $j\ge1$. By \Cref{lem:centering fixed fraction} applied with $\Delta =\Delta _j$, there is
    $n_j' \in \N$ such that $|m_n|/b_n<\varepsilon_j$ for all $n\ge n_j'$.

    By \Cref{lem:tail kill fixed fraction} and \Cref{lem:exp tail fixed fraction} applied with $\Delta =\Delta _j$, there is a projection $p_j$ and $n_j'' \in \N$ such that $\tau(1-p_j)\le \varepsilon_j$, and $v_np_j=0$ for all $n \ge n_j''$ and 
    \begin{align*}
        \limsup_n\normi{(y_n/b_n)p_j}\le\varepsilon_j
    \end{align*}
    
    Now using \eqref{eqn:fixed fraction decomp} and the triangle inequality, for $n\ge\max(n_j',n_j'')$ we have, 
    $ \left\|\frac{S_n(A)}{b_n}p_j\right\|_\infty  \le \left\|\frac{y_n}{b_n}p_j\right\|_\infty+\varepsilon_j+0.$
    Therefore, we have $ \limsup_{n\to\infty}\normi{\frac{S_n(A)}{b_n}p_j} \le 2\varepsilon_j. $ Now fix $\delta>0$ and choose $N\in \N$ with $\sum_{j\ge N}2\varepsilon_j<\delta$. Set $ q_N:=\bigwedge_{j\ge N}p_j \in\CP(\CM). $ Now note that, $\tau(1-q_N)\le\sum_{j\ge N}\tau(1-p_j)<\sum_{j\ge N}2\varepsilon_j<\delta$. Since $q_N\le p_j$ for every $j\ge N$, we have 
    $$ \limsup_{n\to\infty}\normi{\frac{S_n(A)}{b_n}q_N} \le 2\varepsilon_j~~\text{for every }j\ge N. $$
    Note that, $\lim\limits_{j\to \infty} \varepsilon_j=0,$ hence $\lim_{n\to\infty}\normi{\frac{S_n(A)}{b_n}q_N}=0.$ Therefore, $\frac{S_n(A)}{b_n}\ \xrightarrow{\ a.u.\ }\ 0,~~\text{as }n\to\infty,$
\end{proof}

\subsection{Proof of \Cref{mainthm:LIL deterministic weights}(b)}
Fix $\delta\in(0,1]$ and $\Delta>0$. Similar as \Cref{subsec:part(a)}, for $1\le k\le n$ set
$x_k^{(c_n)}:=x_k1_{[0,c_n)}(|x_k|),~ c_n:=\delta b_n/\log n,$
and write $x_k=x_k^{(c_n)}+x_k1_{[c_n,\Delta b_n)}(|x_k|)+x_k1_{[\Delta b_n,\infty)}(|x_k|)$. Then we have, 
\begin{equation}\label{eqn:S_n(A) decompose}
    S_n(A)= y_n+ m_n+w_n+v_n,
\end{equation}
where $y_n:=\sum_ka_{n,k}(x_k^{(c_n)}-\tau(x_k^{(c_n)})), ~m_n:=\sum_ka_{n,k}\tau(x_k^{(c_n)}),~ w_n:=\sum_ka_{n,k}x_k1_{[c_n,\Delta b_n)}(|x_k|)$ and $v_n:=\sum_ka_{n,k}x_k1_{[\Delta b_n,\infty)}(|x_k|).$

If $K=0,$ then $S_n(A)=0$ for every $n$, and the conclusion is immediate. Hence assume $K>0.$

\begin{lem}\label{lem:v_n part}
    For every fixed $\Delta>0$ and every $\varepsilon_1>0$, there is $p_v\in\CP(\CM)$ and $n_0\in N$ such that $\tau(1-p_v)<\varepsilon_1$ and $v_np_v=0$ for all $n\ge n_0.$
\end{lem}
\begin{proof}
    For the proof see \Cref{lem:tail kill fixed fraction}.
\end{proof}

\begin{lem}\label{lem:m_n part}
    Using the definition of $m_n$ from \eqref{eqn:S_n(A) decompose}, we have $\lim\limits_{n\to \infty}\dfrac{m_n}{b_n}=0.$
\end{lem}
\begin{proof}
    Using the definition of $x_k^{(c_n)}$, notice that $\tau(x_k^{(c_n)})=-\tau(x_k1_{[c_n,\infty)}(|x_k|)).$ Therefore,
    $$\frac{\abs{m_n}}{b_n}\le K\cdot \frac{n}{b_n}\int_{[c_n,\infty)}sd\mu(s).$$
    By \Cref{prop:integral s is zero}, it follows that $\lim\limits_{n\to \infty}\frac{m_n}{b_n}=0.$
\end{proof}

\begin{lem}\label{lem:y_n part}
    Let $\theta\in (0,1]$ and $\eta>0.$ Then there exists $\delta=\delta(\theta,\eta)\in(0,1]$ such that, if $c_n=\frac{\delta b_n}{\log n}$ and $y_n$ as in \eqref{eqn:S_n(A) decompose} then the following holds. For every $\varepsilon_0>0,$ there exists $p_y\in \CP(\CM)$ such that $\tau(1-p_y)<\varepsilon_0$ and 
    $$\limsup_{n\to\infty}\frac{\|y_np_y\|_\infty}{b_n}\le(1+\theta)\sqrt{2} K\sqrt{\ell^2+\eta}.$$
\end{lem}
\begin{proof}
    By \Cref{lem:condition verificartion} with $\delta_n=\delta$, it follows that $\|d_{n,k}\|_\infty\le M_n$ and $\sum_kE_{k-1}(d_{n,k}^2)\le D_n^2\cdot \mathbf 1$, where $M_n:2\delta Kb_n/\log n$ and $D_n^2:=K^2n H(c_n)$. Since $0<\delta\le 1$, we have $c_n \le \frac{b_n}{\log n}$. As $H$ is non-decreasing, $H(c_n) \le H\left(\frac{b_n}{\log n}\right).$ Hence, by condition \eqref{item:alpha half condition} of \Cref{mainthm:LIL deterministic weights}, there exists $n_0=n_0(\eta)$ such that, for all $n\ge n_0$, $ H(c_n) \le \dfrac{L_\eta^2b_n^2}{n\log n},$ where we take $L_\eta:=\sqrt{\ell^2+\eta}.$ Consequently, $D_n^2 \le K^2L_\eta^2\dfrac{b_n^2}{\log n}, ~ n\ge n_0.$
    
    Fix $\theta>0,$ and define
    $t_{\theta,\eta} := (1+\theta)\sqrt{2}K L_\eta$ and $ \lambda_n := \frac{\sqrt{2}\log n}{K L_\eta b_n}.$
    Then, notice that $\lambda_n M_n = \dfrac{2\sqrt{2}\delta}{L_\eta}.$ Now choose, $0<\delta \le \min\left\{ 1, 3\theta L_\eta/2\sqrt{2}\right\}.$ Then $\lambda_nM_n\le 3\theta.$ Recall that $y_n=\sum_{k=1}^n d_{n,k}.$ We now apply \Cref{thm: exp ineq JFA}, $ \tau\left(e^{\lambda_n y_n}\right) \le \exp \lrp{\frac{1+\theta}{2} \lambda_n^2D_n^2}.$
    
    Using the bound of $D_n^2,$ for all $n\ge n_0$ we have
    $\frac{1+\theta}{2} \lambda_n^2D_n^2 \le  (1+\theta)\log n.$
    Moreover, $\lambda_nt_{\theta,\eta}b_n  = 2(1+\theta)\log n.$
    Hence, using functional calculus, we have $\tau\left(1_{(t_{\theta,\eta}b_n,\infty)}(y_n)\right)\le n^{-(1+\theta)}.$ 
    The same argument applied to $-y_n$. Hence, $\sum_n\tau(1_{(t_{\theta,\eta}b_n,\infty)}(|y_n|))\le2\sum_nn^{-(1+\theta)}<\infty$. Let $\varepsilon_0>0$. Then, there exists $N=N(\varepsilon_0)$ such that $\sum_{n\ge N} \tau\left( 1_{(t_{\theta,\eta}b_n,\infty)}(|y_n|) \right) <\varepsilon_0.$ Define $p_y := \bigwedge_{n\ge N} 1_{[0,t_{\theta,\eta}b_n]}(|y_n|).$ Then, $\tau(1-p_y)<\varepsilon_0$. Also functional calculus gives 
    $\|y_np_y\|_\infty \le t_{\theta,\eta}b_n, ~~ n\ge N.$ Hence,
    $$\limsup_{n\to\infty}\frac{\|y_np_y\|_\infty}{b_n}\le(1+\theta)\sqrt{2} K\sqrt{\ell^2+\eta}.$$
\end{proof}

The following estimate is inspired by \cite[Lemma 2.3]{cuzick1995strong}. The main idea is to use dyadic blocks $2^{k-1}<n\le 2^k$ and to use the same truncation levels $c_k^*$ and $\Delta b_k^*$ on each block. We then estimate the truncated second moment using $G$ and regular variation. In the critical case $\alpha=1/2$, a direct application of Karamata's theorem is not possible, since the function $th(t)$ has index $-1$. To avoid this endpoint difficulty, we compare $th(t)$ with $(\Delta b_k^*)^{1/2}t^{1/2}h(t)$. The function $t^{1/2}h(t)$ has index $-3/2$, so Karamata's theorem can be applied.

\begin{lem}\label{lem:cuzick's methodology}
    Fix $\beta>0$ and $\Delta>0$. For $k\ge1$, define $b_k^*:=B(2^k), ~ c_k^*:=\beta b_k^*/(k\log 2).$ Then there exists $k_1=k_1(\beta,\Delta)$ such that $c_k^*<\Delta b_k^*,~\forall k\ge k_1.$ For $k\ge k_1$, set
    $$D_k^2 := 2^k\int_{[c_k^*,\,\Delta b_k^*)} s^2\,d\mu(s).$$
    If $\alpha=\frac12$, then there exists a constant $C=C(\beta,\Delta)<\infty$ such that
    $\frac{D_k^2}{(b_k^*)^2} \le Ck^{-1/2}, ~~ k\ge k_1.$
    Consequently, $\displaystyle\sum_{k=1}^{\infty}\left(\dfrac{D_k^2}{(b_k^*)^2}\right)^{M_0}<\infty,$ for every $M_0>2.$
\end{lem}
\begin{proof}
    By \Cref{lem:b_n c_n relation,lem: h properties}, applied with $\delta=\beta$ and $n=2^k,$ we have
    \begin{equation}\label{eqn:h equiv K}
        c_k^*\to \infty ~\text{ and }~c_k^*h(c_k^*)\sim K\left(\frac12,\beta\right)\frac{b_k^*}{2^k}.
    \end{equation}
    Now notice that, by \eqref{eq:Hdef} we have $\int_{[c_k^*,\Delta b_k^*)}s^2\,d\mu(s) = H(\Delta b_k^*)-H(c_k^*).$ Hence, by \Cref{lem:G H integral expressions}(i), note that
    $H(\Delta b_k^*)-H(c_k^*) \le c_k^{*2}G(c_k^*) + \int_{c_k^*}^{\Delta b_k^*}2tG(t)\,dt.$
    Define, $\eta(c):=\sup_{t\ge c}\frac{G(t)}{h(t)}.$ By condition \eqref{item:second condition} of \Cref{mainthm:LIL deterministic weights}, $\lim\limits_{c\to \infty} \eta(c)=0.$ Now, for every $t\ge c_k^*$ we have $G(t)\le\eta(c_k^*)h(t).$ Thus,
    \begin{equation}\label{eqn:int s ineq}
        \int_{[c_k^*,\Delta b_k^*)}s^2\,d\mu(s) \le \eta(c_k^*) \lrp{c_k^{*2}h(c_k^*) + \int_{c_k^*}^{\Delta b_k^*}2th(t)\,dt}.
    \end{equation}
    Since, $\alpha=1/2$, $th(t)\in \mathrm{RV}_{-1}$. Hence, Karamata's theorem can not be applied. We will now use a simple trick to avoid this problem.\\
    For $c_k^*\le t<\Delta b_k^*,$ we have $t\le\Delta b_k^*.$ Hence, $t = t^{1/2}t^{1/2} \le (\Delta b_k^*)^{1/2}t^{1/2}.$ Now observe that,
    \begin{align}
        \int_{c_k^*}^{\Delta b_k^*}2th(t)\,dt
        &\le
        2(\Delta b_k^*)^{1/2}
        \int_{c_k^*}^{\Delta b_k^*}t^{1/2}h(t)\,dt
        \le
        2(\Delta b_k^*)^{1/2}
        \int_{c_k^*}^{\infty}t^{1/2}h(t)\,dt.\label{eqn:2t h ineq}
    \end{align}
    Now take, $\phi(t):=t^{1/2}h(t).$ Since $\alpha=1/2,$ $\phi$ has index $-3/2.$ Note that, $\phi$ has index strictly smaller than $-1$. Thus, by Fact \ref{itm:RV5}, for every $\varepsilon>0$ there exists $c_\varepsilon>0$ such that, for every $c\ge c_\varepsilon,~$ 
    $\int_c^\infty\phi(t)\,dt\le2(1+\varepsilon)c^{3/2}h(c).$
    Since $c_k^*\to\infty$, there exists $k_\varepsilon$ such that $c_k^*\ge c_\varepsilon, ~ k\ge k_\varepsilon.$ Therefore, 
    \begin{equation}\label{eqn: t half h ineq}
        \int_{c_k^*}^{\infty}t^{1/2}h(t)\,dt\le2(1+\varepsilon)(c_k^*)^{3/2}h(c_k^*).
    \end{equation}
    Now using the fact $c_k^*h(c_k^*)\sim K\left(\frac12,\beta\right)\frac{b_k^*}{2^k}$, there exists $C_1>0$ such that for all $k\ge k_\varepsilon$, we have 
    $(c_k^*)^{3/2}h(c_k^*) \le C_1(c_k^*)^{1/2}\frac{b_k^*}{2^k}.$
    Since $c_k^* = \beta b_k^*/(k\log2),$ then for some constant $C_2>0$ we have
    $(c_k^*)^{3/2}h(c_k^*) \le C_2 (b_k^*)^{3/2}/ (2^k k^{1/2}).$ Hence, using \eqref{eqn: t half h ineq}, there exists $C_3>0$ such that$\int_{c_k^*}^{\infty}t^{1/2}h(t)\,dt \le C_3 \frac{(b_k^*)^{3/2}} {2^k k^{1/2}},$ for all $k\ge k_\varepsilon.$ 
    Now, combining with \eqref{eqn:2t h ineq}, we have
    \begin{equation}\label{eqn:2t h(t) ineq}
        \int_{c_k^*}^{\Delta b_k^*}2th(t)\,dt \le C_4\cdot  \frac{(b_k^*)^2} {2^k k^{1/2}},
    \end{equation}
    for some constant $C_4>0.$ Now from \eqref{eqn:h equiv K}, there exists $k_0\in \N$ such that for all $k\ge k_0$, we have $c_k^{*2}h(c_k^*)= c_k^*\lrp{c_k^*h(c_k^*)}\le C_5 \frac{(b_k^*)^2}{2^k k}$, for some constant $C_5>0.$ Thus using \eqref{eqn:2t h(t) ineq}, there exists $k_0<k_1\in \N$ and $C_6>0$ such that
    $$c_k^{*2}h(c_k^*) + \int_{c_k^*}^{\Delta b_k^*}2th(t)\,dt \le C_6 \frac{(b_k^*)^2} {2^k k^{1/2}}.$$
    Since $\eta(c_k^*)\to0$, there exists $k_2$ such that, $\eta(c_k^*)\le1,~~ k\ge k_2.$ Combining this with \eqref{eqn:int s ineq}, we obtain, $\int_{[c_k^*,\Delta b_k^*)}s^2\,d\mu(s) \le C_6 \cdot \frac{(b_k^*)^2} {2^k k^{1/2}}$, for $k\ge k_3:= \max\{k_1,k_2\}.$ 
    Since $D_k^2 = 2^k \int_{[c_k^*,\Delta b_k^*)}s^2\,d\mu(s)$, for all $k\ge k_3,$ we have $D_k^2/(b_k^*)^2\le C_6 k^{-1/2}$. Finally, for $M_0>2$ we have $\lrp{D_k^2/(b_k^*)^2}^{M_0}\le C k^{-M_0/2},$ where $C=C_6^{M_0}.$ Since $M_0>2,$ $\displaystyle\sum_{k=1}^{\infty}\left(\dfrac{D_k^2}{(b_k^*)^2}\right)^{M_0}<\infty.$    
\end{proof}

The following lemma provides a non-commutative counterpart of \cite[Lemma 2.2]{cuzick1995strong} and will play a key role in the proof of the intermediate truncation term. 

\begin{lem}\label{lem:NC_Cuzick_2.2}
    Let $ d_{k,i}\in\CM_{k,i},~ 1\le i\le 2^k, $ be self-adjoint elements of $\CM$ satisfying $ E_{k,i-1}(d_{k,i})=0, ~ \|d_{k,i}\|_\infty\le b_k ~ (1\le i\le 2^k), $ where $(b_k)_{k\ge1}$ is a sequence of positive reals with $b_k\uparrow\infty$. Put $ S_k:=\sum_{i=1}^{2^k}d_{k,i}. $ Assume that there is a sequence $(D_k)_{k\ge1}$ of nonnegative reals such that $\sum_{i=1}^{2^k}E_{k,i-1}\lrp{d_{k,i}^2} \le D_k^2\cdot 1$ and $\sum_{k=1}^{\infty}\left(D_k^2/b_k^2\right)^{M}<\infty,$ for some constant $M>0$. 
    Then, $\limsup\limits_{k\to\infty}\frac{S_k}{b_k}\ \underset{\mathrm{a.u.}}{\le}\ M.$
\end{lem}
\begin{proof}
    Let $r_k:=D_k^2/b_k^2,~ k\ge1.$ By hypothesis, $\sum_k r_k^M<\infty$. Therefore,  $r_k\to0$ as $k\to\infty$. Fix $\varepsilon,\delta>0$ and set $a:=M+\delta$. Fix $k$ with $r_k>0$. Define $x_{k,i}:=\sum_{j=1}^{i}d_{k,j}$ for $0\le i\le 2^k$, so that $x_{k,0}=0$ and $x_{k,2^k}=S_k$. $(x_{k,i})_{i=0}^{2^k}$ is a self-adjoint martingale with respect to $(\CM_{k,i},E_{k,i})_{i=0}^{2^k}$, with differences $d_{k,i}$ satisfying $\|d_{k,i}\|_\infty\le b_k$ and, $\sum_iE_{k,i-1}(d_{k,i}^2)\le D_k^2\cdot1$. Now, applying \Cref{thm: exp ineq JFA}, for every $\lambda>0$,
    \begin{equation}\label{eqn:JFA ineq here}
        \tau\lrp{e^{\lambda S_k}} \le \exp\lrp{r_k(e^{\lambda b_k}-\lambda b_k-1)}.
    \end{equation}
    Define $q_k:=1_{(ab_k,\infty)}\lrp{|S_k|}$. Hence using \cref{eqn:JFA ineq here},
    $\tau\lrp{1_{(ab_k,\infty)}(S_k)}\ \le \exp\lrp{-\lambda ab_k+r_k(e^{\lambda b_k}-\lambda b_k-1)}.$
    Writing $t:=\lambda b_k>0$, 
    \begin{equation} \label{eqn: q_k f_k ineq}
        \tau(q_k) \le 2\exp(f_k(t)),~~ f_k(t):=-at+r_k(e^t-t-1),~ t>0.
    \end{equation}
    $f_k$ attains its minimum at $t_k^*=\log\lrp{1+\frac{a}{r_k}}>0. $ Fix $\theta\in(0,\tfrac12)$ such that $a(1-2\theta)>M$. Since $r_k\to0$, there is $k_1$ such that $r_k<1$ for all $k\ge k_1$. Hence 
    $t_k^*\ \ge \log a+\log\frac1{r_k},~~\text{ for } k\ge k_1.$ 
    Since $\log(1/r_k)\to\infty$, there is $k_2\ge k_1$ such that $\log a\ge -\theta\log(1/r_k)$ for all $k\ge k_2$, giving 
    \begin{equation}\label{eqn:tk star ineq}
        t_k^* \ge (1-\theta)\log\frac{1}{r_k},~~\text{ for } k\ge k_2. 
    \end{equation}
    Hence by \eqref{eqn:tk star ineq} and $t_k^*=\log\lrp{1+\frac{a}{r_k}}$, we have
    $f_k(t_k^*)  \le a-a(1-\theta)\log\frac1{r_k},~~\text{ for } k\ge k_2.$
    Since $\log(1/r_k)\to\infty$, choose $k_3\ge k_2$ so that $a\le\theta a\log(1/r_k)$ for $k\ge k_3$. Then we have, $f_k(t_k^*)\le -M\log\frac1{r_k}$ for $k\ge k_3$. Hence, using \eqref{eqn: q_k f_k ineq}, we get $\tau(q_k)\le 2\,r_k^{M}$, for $k\ge k_3$.\\
    Now by hypothesis and above argument we have, $\sum_{k\ge k_3}\tau(q_k)\le 2\sum_{k\ge k_3}r_k^M<\infty$. Choose $N\ge k_3$ such that $ \sum_{k\ge N}\tau(q_k)<\varepsilon. $ Define $ p:=\bigwedge_{k\ge N}(1-q_k)\ \in\ \CP(\CM). $ Fix $k\ge N$. Then $S_kp\ =\ S_k(1-q_k)p$ and  $\|S_k(1-q_k)\|_\infty\le ab_k$. Therefore, $\left\|\frac{S_k}{b_k}p\right\|_\infty\le a = M+\delta,$ for all $k\ge N$.
    Hence, $\limsup_{k\to\infty}\frac{S_k}{b_k}\ \underset{\mathrm{a.u.}}{\le}\ M.$
\end{proof}

The dyadic decomposition used in the following proof is inspired by \cite[Lemma 2.3]{cuzick1995strong}. In particular, the main idea is to divide the integers into dyadic blocks $2^{k-1}<n\le 2^k$ and use the same interval on each block.

\begin{lem}\label{lem:w_n part}
    Let $\theta,\eta$ be as in \Cref{lem:y_n part}, and let $\delta$ be the corresponding constant. For every $\rho>0$, there exists $\Delta=\Delta(\rho)>0$ such that, for every $\varepsilon_1>0$, there is a projection $p_w\in\CP(\CM)$ satisfying $ \tau(1-p_w)<\varepsilon_1$ and $$\limsup_{n\to\infty} \frac{\|p_w w_np_w\|_\infty}{b_n} \le \rho.$$
\end{lem}
\begin{proof}
    Recall that $b_k^*= B(2^k)$ for all $k \in \N$. By \Cref{lem: property1 of RV func}, applied to $B\in\mathrm{RV}_{1/2}$ with the fixed ratio $\lambda=1/2$, there exists $C_B\ge1$ and $k_0\in\mathbb N$ such that
    $C_B^{-1}b_k^* \le b_n \le b_k^* ~\text{ and }~  2^{k-1}<n\le2^k, ~~\text{for } k\ge k_0.$
    Following \Cref{lem:cuzick's methodology}, set $\beta:=\delta/C_B, ~ c_k^*:=\beta b_k^*/(k\log2).$ Fix $k\ge k_0$ and let $2^{k-1}<n\le2^k.$ Since $\log n\le k\log2$, we have 
    $$c_n =\frac{\delta b_n}{\log n} \ge \frac{\delta C_B^{-1}b_k^*}{k\log2} =c_k^*.$$
    Also, $\Delta b_n\le\Delta b_k^*.$ Hence, $[c_n,\Delta b_n) \subseteq [c_k^*,\Delta b_k^*).$ For $1\le i\le2^k$, define
    $$z_{i,k} := |x_i| 1_{[c_k^*,\Delta b_k^*)}(|x_i|) - \tau\lrp{ |x_i|1_{[c_k^*,\Delta b_k^*)}(|x_i|)} \mathbf 1.$$
    Set, $u_k:=\sum_{i=1}^{2^k}z_{i,k}$ and $h_k := \sum_{i=1}^{2^k} |x_i|1_{[c_k^*,\Delta b_k^*)}(|x_i|).$ Since $(x_i)$ are i.i.d, writing $a_k:= 2^k \tau\lrp{|x_1|1_{[c_k^*,\Delta b_k^*)}(|x_1|)},$ we have $h_k=u_k+ a_k\cdot \mathbf 1.$ By \Cref{lem:y_k d_{n,k} introduce}, we have $E_{k,i-1}(z_{i,k})=0,$ $\|z_{i,k}\|_\infty \le 2\Delta b_k^*,$ and $\sum_{i=1}^{2^k} E_{k,i-1}(z_{i,k}^2) \le D_k^2\cdot \mathbf 1,$ where $D_k^2$ is from \Cref{lem:cuzick's methodology}. By \Cref{lem:cuzick's methodology}, with the above choices of $\beta$ and $\Delta$, we have $\sum\limits_{k=1}^{\infty} \lrp{\dfrac{D_k^2}{(2\Delta b_k^*)^2} }^3 <\infty.$ Therefore, \Cref{lem:NC_Cuzick_2.2} applies to the row $(z_{i,k})_{1\le i\le2^k}$ with bound $2\Delta b_k^*$. Thus, for every $\varepsilon_1>0$, there exists a projection $p_w\in\CP(\CM)$ such that $\tau(1-p_w)<\varepsilon_1$ and $\limsup_{k\to\infty} \|u_kp_w\|_\infty/2\Delta b_k^* \le3.$ Hence,
    \begin{equation}\label{eqn: u_k p_w norm ineq}
        \limsup_{k\to\infty}\frac{\|u_kp_w\|_\infty}{b_k^*}\le6\Delta.
    \end{equation}
    Next, notice that
    $$0\le \frac{a_k}{b_k^*} = \frac{2^k}{b_k^*} \int_{[c_k^*,\Delta b_k^*)} s\, d\mu(s) \le \frac{2^k}{b_k^*} \int_{[c_k^*,\infty)} s\, d\mu(s).$$
    Now \Cref{prop:integral s is zero} with $n=2^k$ and $\delta=\beta$ gives $a_k/b_k^*\longrightarrow0$ as $k\to \infty.$ Since $a_k$ is a scalar, $\|a_kp_w\|_\infty=a_k.$ Therefore, from $h_k=u_k+a_k1$ and \eqref{eqn: u_k p_w norm ineq}, we have
    \begin{equation}\label{eqn: h_k ineq}
        \limsup_{k\to\infty}\frac{\|h_kp_w\|_\infty}{b_k^*}\le\limsup_{k\to\infty}\frac{\|u_kp_w\|_\infty}{b_k^*}+\limsup_{k\to\infty}\frac{a_k}{b_k^*}\le6\Delta.
    \end{equation}
    Let $2^{k-1}<n\le2^k.$ Since, $[c_n,\Delta b_n) \subseteq [c_k^*,\Delta b_k^*)$ and $K=\sup_{n,k}\abs{a_{n,k}}$ we have, 
    \begin{align}
        -Kh_k\le w_n\le Kh_k \implies \| p_ww_np_w\|_\infty \le K\|p_wh_k p_w\|_\infty \le K\normi{h_k p_w}\label{eqn: norm p_w ineq}
    \end{align}
    Since, $b_n\ge C_B^{-1}b_k^*,$ we have $b_n^{-1}\le C_B/b_k^*.$ Combining this with \eqref{eqn: norm p_w ineq} and \eqref{eqn: h_k ineq} we have \\
    $\limsup\limits_{n\to\infty}\|p_ww_np_w\|_\infty/b_n \le 6KC_B\Delta.$
    Now, choose $ \Delta:=\frac{\rho}{6KC_B}.$ Then, $\limsup\limits_{n\to\infty}\dfrac{\|p_ww_np_w\|_\infty}{b_n}\le\rho.$ Together with $\tau(1-p_w)<\varepsilon_1$, this proves the lemma.
\end{proof}

\begin{proof}[Proof of \Cref{mainthm:LIL deterministic weights}(b)]
    Fix $\varepsilon_0,\varepsilon>0$. Choose $\theta\in(0,1]$, $\eta>0$, $\rho>0$ so small that $(1+\theta)\sqrt2K\sqrt{\ell^2+\eta}+\rho<\sqrt2K\ell+\varepsilon.$
    Let $\delta$ and $\Delta=\Delta(\rho)$ be as given by \Cref{lem:y_n part} and \Cref{lem:w_n part}, respectively. Now notice that, using \Cref{lem:y_n part}, there exists a projection $p_y\in\CP(\CM)$ and $N_y\in \N$ such that $\tau(1-p_y)<\frac{\varepsilon_0}{3}$ and 
    \begin{equation}\label{eqn:y_n part}
        \normi{y_n p_y}\le (1+\theta)\sqrt{2}K\sqrt{\ell^2+\eta}\,b_n,~~\forall n\ge N_y.
    \end{equation}
    By \Cref{lem:w_n part}, there exist $p_w\in\CP(\CM)$ and $N_w\ge1$  such that $\tau(1-p_w)<\frac{\varepsilon_0}{3}$ and 
    \begin{equation}\label{eqn: w_n p_w part}
        \frac{\|p_ww_np_w\|_\infty}{b_n}\le\rho+\varepsilon,~~\forall n\ge N_w.
    \end{equation}
    By \Cref{lem:v_n part}, there exist $p_v\in\CP(\CM)$ and $N_v\ge1$ such that $\tau(1-p_v)<\frac{\varepsilon_0}{3}$ and 
    \begin{equation}\label{eqn:v_n part}
        p_vv_np_v=0, ~~\forall n\ge N_v.
    \end{equation}
    Now, define $p:=p_y\wedge p_w\wedge p_v.$ Then, $\tau(1-p)<\varepsilon_0.$ Set, $N:=\max\{N_y,N_v,N_w\}.$ Therefore, using \eqref{eqn:S_n(A) decompose} and \eqref{eqn:v_n part}, for every $n\ge N$, $\|pS_n(A)p\|_\infty\le\|y_np_y\|_\infty+|m_n|+\|p_ww_np_w\|_\infty.$    Now using \Cref{eqn:y_n part,eqn: w_n p_w part,lem:m_n part}, we have for every $n\ge N,$ we have
    $$\frac{\|pS_n(A)p\|_\infty}{b_n}\le(1+\theta)\sqrt{2}K\sqrt{\ell^2+\eta}+\frac{|m_n|}{b_n}+\frac{\|p_ww_np_w\|_\infty}{b_n}.$$
    Taking limsup both side, and using \Cref{lem:m_n part,eqn: w_n p_w part}, we have
    $\limsup_{n\to\infty}\frac{\|pS_n(A)p\|_\infty}{b_n}\le \sqrt{2}K\ell+\varepsilon.$
    Since $\varepsilon,\varepsilon_0>0$ are arbitrary, we have
    $$\limsup_{n\to\infty}\frac{S_n(A)}{b_n}\underset{\mathrm{b.a.u.}}{\le}\sqrt{2}K\ell.$$
\end{proof}

\subsection*{Acknowledgement}

DS acknowledges the Institute of Mathematics of the Polish Academy of Sciences, for providing the institutional support during the completion of this work. SP gratefully acknowledges the Indian Statistical Institute for financial support through its doctoral fellowship. He also wishes to thank his PhD supervisor, Prof. Issan Patri, for his valuable guidance and support.


\end{document}